%%%%%%%%%%%%%%%%%%%%%%%%%%%%%%%%%%%%%%%%%%%%%%%%%%%%%%%%%%
%%
%%     cgps3 su viscous Cahn-Hilliard 
%%
%%%%%%%%%%%%%%%%%%%%%%%%%%%%%%%%%%%%%%%%%%%%%%%%%%%%%%%%%%

\def\LaTeX{%
  \let\Begin\begin
  \let\End\end
  \let\salta\relax
  \let\finqui\relax
  \let\futuro\relax}

\def\UK{\def\our{our}\let\sz s}
\def\USA{\def\our{or}\let\sz z}

%%%%%%%%%%%%%%%%%%%%%%%%%%%%%%%%%

% scegliere fra \TeX e \LaTeX  e fra  \UK oppure \USA

%\TeX
\LaTeX

%\UK
\USA

%%%%%%%%%%%%%%%%%%%%%%%%%%%%%%%%%
%% page layout
%%%%%%%%%%%%%%%%%%%%%%%%%%%%%%%%%

\salta

\documentclass[twoside,12pt]{article}
\setlength{\textheight}{24cm}
\setlength{\textwidth}{16cm}
\setlength{\oddsidemargin}{2mm}
\setlength{\evensidemargin}{2mm}
\setlength{\topmargin}{-15mm}
\parskip2mm

%%%%%%%%%%%%%%%%%%%%%%%%%%%%%%%%%
%% packages
%%%%%%%%%%%%%%%%%%%%%%%%%%%%%%%%%

\usepackage{color}
\usepackage{amsmath}
\usepackage{amsthm}
\usepackage{amssymb}
\usepackage[mathcal]{euscript}

%\usepackage[notref,notcite]{showkeys}
%\usepackage{showkeys}

%\newcommand{\red}{\marginpar{{$\leftarrow$}}\color{red}}
%\newcommand{\green}{\marginpar{{$\leftarrow$}}\color{green}}
%\def\green{\color{green}}
%\def\red{\color{red}}
%\def\Green #1{{\green #1}}
%\def\red{}
%\def\green{}
%\def\Green #1{}

%%%%%%%%%%%%%%%%%%%%%%%%%%%%%%%%%
%% bibliographystyle
%%%%%%%%%%%%%%%%%%%%%%%%%%%%%%%%%

\bibliographystyle{plain}

%%%%%%%%%%%%%%%%%%%%%%%%%%%%%%%%%
%% environments
%%%%%%%%%%%%%%%%%%%%%%%%%%%%%%%%%

%

\finqui

\def\Beq{\Begin{equation}}
\def\Eeq{\End{equation}}
\def\Bsist{\Begin{eqnarray}}
\def\Esist{\End{eqnarray}}

\def\Bthm{\Begin{theorem}}
\def\Ethm{\End{theorem}}
\def\Blem{\Begin{lemma}}
\def\Elem{\End{lemma}}

\def\Brem{\Begin{remark}\rm}
\def\Erem{\End{remark}}

\let\non\nonumber

%%%%%%%%%%%%%%%%%%%%%%%%%%%%%%%%%
%% macros
%%%%%%%%%%%%%%%%%%%%%%%%%%%%%%%%%

% macro salvate

\newcommand\QED{\hfill $\square$}

% sottosezioni non numerate

\def\step #1 \par{\medskip\noindent{\bf #1.}\quad}

% abbreviazioni di parole

\def\holder{H\"older}
\def\aand{\quad\hbox{and}\quad}
\def\loti{long-time}
\def\Loti{Long-time}

\def\lhs{left-hand side}
\def\rhs{right-hand side}
\def\sfw{straightforward}
\def\wk{well-known}
\def\wepo{well-posed}
\def\Wepo{Well-posed}

% versioni inglesi (UK) o americane (USA)

\def\characteriz{characteri\sz}

\def\organiz{organi\sz}

\def\regulariz{regulari\sz}

\def\bhv{behavi\our}

% bold, cal e mathop

\def\multibold #1{\def\arg{#1}%
  \ifx\arg\pto \let\next\relax
  \else
  \def\next{\expandafter
    \def\csname #1#1#1\endcsname{{\bf #1}}%
    \multibold}%
  \fi \next}

\def\pto{.}

\def\multical #1{\def\arg{#1}%
  \ifx\arg\pto \let\next\relax
  \else
  \def\next{\expandafter
    \def\csname cal#1\endcsname{{\cal #1}}%
    \multical}%
  \fi \next}

% operatori

\def\multimathop #1 {\def\arg{#1}%
  \ifx\arg\pto \let\next\relax
  \else
  \def\next{\expandafter
    \def\csname #1\endcsname{\mathop{\rm #1}\nolimits}%
    \multimathop}%
  \fi \next}

\multibold
qwertyuiopasdfghjklzxcvbnmQWERTYUIOPASDFGHJKLZXCVBNM.

\multical
QWERTYUIOPASDFGHJKLZXCVBNM.

\multimathop
dist div dom meas sign supp .

% accorpamenti di formule citate:
% uso  \accorpa {prima}{seconda}
%      \Accorpa\cs prima seconda (con il comodo blank anche dopo)
% NB: \Accorpa definisce \cs come l'accorpamento delle due citazioni
% e scrive sul file.log

\def\accorpa #1#2{\eqref{#1}--\eqref{#2}}
\def\Accorpa #1#2 #3 {\gdef #1{\eqref{#2}--\eqref{#3}}%
  \wlog{}\wlog{\string #1 -> #2 - #3}\wlog{}}

% macro comode

\def\graffe #1{\mathopen\{#1\mathclose\}}

\def\<#1>{\mathopen\langle #1\mathclose\rangle}
\def\norma #1{\mathopen \| #1\mathclose \|}

\def\iot {\int_0^t}
\def\ioT {\int_0^T}
\def\iO{\int_\Omega}
\def\intQt{\iot\!\!\iO}
\def\intQ{\ioT\!\!\iO}

\def\dt{\partial_t}
\def\dn{\partial_\nu}

\def\cpto{\,\cdot\,}

\def\checkmmode #1{\relax\ifmmode\hbox{#1}\else{#1}\fi}
\def\aeO{\checkmmode{a.e.\ in~$\Omega$}}
\def\aeQ{\checkmmode{a.e.\ in~$Q$}}

\def\aaQ{\checkmmode{for a.a.~$(x,t)\in Q$}}
\def\aat{\checkmmode{for a.a.~$t\in(0,T)$}}

% insiemi numerici

\def\erre{{\mathbb{R}}}

% spazi di funzioni a valori vettoriali su [0,T], [0,t], [0,s], [0,+\infty), [\delta,T]

% Come ricordare: in generale i simboli L H W  C da soli per gli spazi su (0,T)
% gli stessi raddoppiati per (0,+\infty)
% aggiunta di t o s al simbolo per (0,t) e (0,s)
% aggiunta di d al simbolo semplice o doppio per intervalli (\delta,T) e (\delta,+\infty)
% il simbolo C e i suoi derivati mettono le quadre anziche' le tonde

% Esempi   \L2V   \L\infty\Vp   \W{1,1}H   \C0H   \LL2V   \CC0\Vp   \Ld2V  \CCdH

\def\genspazio #1#2#3#4#5{#1^{#2}(#5,#4;#3)}
\def\spazio #1#2#3{\genspazio {#1}{#2}{#3}T0}

\def\L {\spazio L}
\def\H {\spazio H}
\def\W {\spazio W}

\def\C #1#2{C^{#1}([0,T];#2)}

% spazi di funzioni su \Omega e \Gamma

\def\Lx #1{L^{#1}(\Omega)}
\def\Hx #1{H^{#1}(\Omega)}

\def\Cx #1{C^{#1}(\overline\Omega)}

\def\Luno{\Lx 1}
\def\Ldue{\Lx 2}
\def\Linfty{\Lx\infty}
\def\Lq{\Lx4}
\def\Huno{\Hx 1}
\def\Hdue{\Hx 2}

% spazi di funzioni su Q

\def\LQ #1{L^{#1}(Q)}

% lettere greche

\let\theta\vartheta
\let\eps\varepsilon

\let\TeXchi\chi                         % new \chi, exactly on the baseline
\newbox\chibox
\setbox0 \hbox{\mathsurround0pt $\TeXchi$}
\setbox\chibox \hbox{\raise\dp0 \box 0 }
\def\chi{\copy\chibox}

% abbreviazioni specifiche del lavoro

\def\muz{\mu_0}
\def\rhoz{\rho_0}
\def\uz{u_0}
\def\muzs{\mu_0^*}
\def\rhomin{\rho_*}

\def\normaV #1{\norma{#1}_V}
\def\normaH #1{\norma{#1}_H}
\def\normaW #1{\norma{#1}_W}
\def\normaVp #1{\norma{#1}_{\Vp}}

\def\T{\calT_\tau}

\def\mut{\mu_\tau}
\def\rhot{\rho_\tau}

\def\rhonl{\rho_n^\lambda}

\def\fl{f_\lambda}
\def\ful{f_{1,\lambda}}
\def\fdl{f_{2,\lambda}}

\def\mun{\mu_n}
\def\rhon{\rho_n}
\def\munmu{\mu_{n-1}}

\def\iotmt{\int_0^{t-\tau}}

\def\chik{\chi_k}
\def\mumkp{(\mu-k)^+}
\def\rhomrm{(\rho-\rmin)^-}
\def\rmin{r_*}
\def\rmax{r^*}

\def\Norma #1{\mathopen{|\!|\!|}#1\mathclose{|\!|\!|}}
\def\Ldlq{\L2\Lq}

\def\mui{\mu_\infty}
\def\rhoi{\rho_\infty}
\def\phii{\phi_\infty}
\def\mus{\mu_s}
\def\rhos{\rho_s}
\def\muo{\mu_\omega}
\def\rhoo{\rho_\omega}

\def\itnt{\int_{t_n}^{t_n+t}}
\def\intQtnt{\itnt\!\!\!\iO}

\def\Vp{V^*}

% scelte che si possono cambiare

\def\Dzero{D_0}
\def\Duno{D_1}
\def\Ddue{D_2}

%%%%%%%%%%%%%%%%%%%
%new definitions introduced by PPG
%%%%%%%%%%%%%%%%%%%
%\usepackage{greekbf}
%\def\chi{{\Greekmath 011F}}%
\DeclareMathAlphabet{\mathbf}{OT1}{cmr}{bx}{it}

\newcommand{\hb}{\mathbf{h}}
\newcommand{\nb}{\mathbf{n}}
\newcommand{\Mb}{\mathbf{M}}
\newcommand{\Hb}{\mathbf{H}}
 \font\mba=cmmib10 scaled
\magstephalf

\def\csib{\hbox{\mba {\char 24}}}
%%%%%%%%%%%%%%%%%%%%%%%%%%%%%%
\Begin{document}
%%%%%%%%%%%%%%%%%%%%%%%%%%%%%%%%%

%%%%%%%%%%%%%%%%%%%%%%%%%%%%%%%%%
%% front page
%%%%%%%%%%%%%%%%%%%%%%%%%%%%%%%%%

\title{{\bf \Wepo ness and \loti\ \bhv\
  for a {nonstandard} viscous Cahn-Hilliard
  {system}}\footnote{{\bf Acknowledgments.}\quad\rm
The authors gratefully acknowledge {the} financial
support {of} the MIUR-PRIN Grant 2008ZKHAHN \emph{``Phase transitions, hysteresis 
and multiscaling''}, the EU Marie Curie Research Training Network
MULTIMAT \emph{``Multi-scale Modeling and Characterization
for Phase Transformations in Advanced Materials''},
the DFG Research Center \emph{Matheon} in Berlin, 
and the IMATI of CNR in Pavia.}}
\author{}
\date{}
\maketitle
\begin{center}
\vskip-2.2cm
{\large\bf Pierluigi Colli$^{(1)}$}\\
{\normalsize e-mail: {\tt pierluigi.colli@unipv.it}}\\[.4cm]
{\large\bf Gianni Gilardi$^{(1)}$}\\
{\normalsize e-mail: {\tt gianni.gilardi@unipv.it}}\\[.4cm]
{\large\bf Paolo Podio-Guidugli$^{(2)}$}\\
{\normalsize e-mail: {\tt ppg@uniroma2.it}}\\[.4cm]
{\large\bf J\"urgen Sprekels$^{(3)}$}\\
{\normalsize e-mail: {\tt sprekels@wias-berlin.de}}\\[.6cm]
$^{(1)}$
{\small Dipartimento di Matematica ``F. Casorati'', Universit\`a di Pavia}\\
{\small via Ferrata 1, 27100 Pavia, Italy}\\[.2cm]
$^{(2)}$
{\small Dipartimento di Ingegneria Civile, Universit\`a di Roma ``Tor Vergata''}\\
{\small via del Politecnico 1, 00133 Roma, Italy}\\[.2cm]
$^{(3)}$
{\small Weierstra\ss-Institut f\"ur Angewandte Analysis und Stochastik}\\
{\small Mohrenstra\ss e\ 39, 10117 Berlin, Germany}\\[.8cm]
\end{center}

\Begin{abstract}
We study a diffusion model of phase field type, consisting  {of}
a system of two partial differential equations {encoding} the balance{s}
of microforces and microenergy; the two unknowns are the order parameter
and the chemical potential. By a careful development of uniform estimates and 
the deduction of {certain} useful boundedness properties, we
prove existence and uniqueness of a global-in-time
smooth solution to the {associated} initial/boundary-value problem; 
moreover, we give a description of the relative $\omega$-limit set.

\vskip3mm

\End{abstract}

\salta

\pagestyle{myheadings}
\newcommand\testopari{\sc Colli \ --- \ Gilardi \ --- \ Podio-Guidugli \ --- \ Sprekels}
\newcommand\testodispari{\sc \Wepo ness and \loti\ for a {nonstandard diffusion} model}
\markboth{\testodispari}{\testopari}

%\markright{}

\finqui

%%%%%%%%%%%%%%%%%%%%%%%%%%%%%%%%%
%% very beginning
%%%%%%%%%%%%%%%%%%%%%%%%%%%%%%%%%

%\section{Introduction}
\section{Problem setting}
\label{Intro}
\setcounter{equation}{0}
The Cahn-Hilliard {system}:
\Beq\label{CH}
 \dt\rho - \kappa \Delta \mu =0 \ , \qquad  \mu= - \Delta\rho + f'(\rho)
 {,}  
 \Eeq
{describes diffusion-driven}
phase{-segregation} processes in a two-phase material body. 
{Here} 
$\rho$, with $\rho(x,t)\in[0,1]$, is an \emph{order parameter} field 
interpreted as the scaled volumetric density of one of the two phases, 
$\kappa>0$~is a \emph{mobility} coefficient, and $\mu$  {is} the \emph{chemical potential}{;}
%{that is related to $\rho$ via the second equation in} Here, 
$f'$ stands for the derivative of a double-well potential $f$. {Customarily, the two equations \eqref{CH} are combined so as to obtain the \emph{Cahn-Hilliard equation}:
\Beq\label{CHe}
\dt\rho = \kappa \Delta (- \Delta\rho + f'(\rho)),
\Eeq
 a nonlinear high-order parabolic PDE for the order parameter that has been studied estensively. With this procedure -- we note for later reference -- the chemical potential is left in the background; in particular, there is no need to take an a priori decision about its sign.}
% {we assume 
%that $f'$ is significantly confined in~$(0,1)$ and singular at the endpoints}. 

{To achieve their generalization of \eqref{CHe}, Fried \& Gurtin and Gurtin \cite{FG,Gurtin} propose: (i) to regard the second of \eqref{CH} as} a \textit{balance of microforces:}
\Beq\label{balance}
\div\csib+\pi+\gamma=0{,}
\Eeq
{where the distance microforce per unit volume is split into} an internal part $\pi$ and an 
external part $\gamma$, and the contact microforce per unit area of a surface oriented by its normal $\nb$ is measured by $\csib\cdot\nb$ in terms of the \emph{microstress} vector  $\csib$;{\footnote{{In \cite{Fremond},  the microforce balance
is stated under the form of a principle of virtual powers 
for microscopic motions.}} (ii) to interpret the first equation as a}  
\textit{balance law for the order parameter}{:}
\Beq\label{balorpam}
\partial_t\rho = - \div \hb + \sigma{,}
\Eeq
where the pair $(\hb ,  \sigma)$ is the \textit{inflow} of $\rho$; {(iii) to restrict the admissible constitutive choices for $\pi,\csib, \hb$, and the \emph{free energy density} $\psi$, to those consistent in the sense of Coleman \& Noll \cite{CN} with an \emph{ad hoc} version of the Second Law of continuum thermodynamics,  namely a postulated ``dissipation inequality that accomodates diffusion'':
\begin{equation}\label{dissipation}
\partial_t\psi +(\pi-\mu)\partial_t\rho-\csib\cdot\nabla(\partial_t\rho)+\hb\cdot\nabla\mu\leq 0
\end{equation}
(cf. eq. (3.6) in \cite{Gurtin}). Within this framework, the following set of constitutive prescriptions is shown acceptable:
\Beq\label{costi}
\begin{aligned}
\psi&=\widehat\psi(\rho,\nabla\rho),\\ \widehat\pi(\rho,\nabla\rho,\mu)&=\mu- \partial_\rho \widehat\psi(\rho,\nabla\rho), \\ \widehat\csib(\rho,\nabla\rho)&=\partial_{\nabla\rho} \widehat\psi(\rho,\nabla\rho),
\end{aligned}
\Eeq
 together with 
 \Beq\label{acca}
\hb = - \Mb\nabla \mu , \quad \hbox{with } \ \Mb=\widehat\Mb(\rho,\nabla\rho,\mu, \nabla\mu);
\Eeq
moreover, it is shown that the tensor-valued \emph{mobility mapping} $\Mb$ must satisfy the inequality:
\Beq
\nabla \mu\cdot \widehat\Mb(\rho,\nabla\rho,\mu, \nabla\mu) \nabla\mu \geq 0 . \non
\Eeq
It follows from \eqref{balance}, \eqref{balorpam}, \eqref{costi}, and $\eqref{acca}_1$ that:
\[
\dt\rho = \div\left(\Mb\nabla\left(\partial_\rho \widehat\psi(\rho,\nabla\rho)-\div\big(\partial_{\nabla\rho} \widehat\psi(\rho,\nabla\rho)\big)-\gamma\right)\right)+\sigma
\]
(cf. eq. (3.17) in \cite{Gurtin}); in particular, the Cahn-Hilliard equation \eqref{CHe} is arrived at by taking:
\Bsist
&&\widehat\psi(\rho,\nabla\rho)= f(\rho)+\frac{1}{2}|\nabla\rho|^2,\qquad \Mb=\kappa \mathbf{1},
\label{constitutive}
\Esist
and both the external distance microforce $\gamma$ and the order-parameter source term $\sigma$ identically null.
}

One of us proposed in \cite{Podio} a modified version of Fried \& Gurtin's derivation{, in which} {their step (i) is retained, but the order-parameter balance \eqref{balorpam} and the dissipation inequality \eqref{dissipation} are both dropped and replaced, respectively, by}
 the \emph{microenergy balance}
\Beq\label{energy}
\partial_t\varepsilon=e+w,\quad e:=-\div{\overline \hb}+{\overline \sigma},\quad w:=-\pi\,\partial_t\rho+\csib\cdot\nabla(\partial_t\rho)
\Eeq
and the \emph{microentropy imbalance}
\Beq\label{entropy}
\partial_t\eta\geq -\div\hb+\sigma,\quad \hb:=\mu{\overline \hb},\quad \sigma:=\mu\,{\overline \sigma}.
\Eeq
The salient new feature of this approach to phase-segregation modeling is that the \emph{microentropy inflow} $(\hb,\sigma)$ is deemed proportional to the \emph{microenergy inflow} $({\overline \hb},{\overline\sigma})$ through the \emph{chemical potential} $\mu$, a {positive} field; consistently, the free energy is defined~to~be
\Beq\label{freeenergy}
\psi:=\varepsilon-\mu^{-1}\eta,
\Eeq
with chemical potential playing the same role as \emph{coldness} in the deduction of the heat equation.\footnote{{\color{black}As much as absolute temperature is a macroscopic measure of microscopic \emph{agitation}, its inverse - the coldness - measures microscopic \emph{quiet}; likewise, as argued in \cite{Podio},  chemical potential can be seen as a macroscopic measure of microscopic \emph{organization}.}}  Combination of (\ref{energy})-(\ref{freeenergy}) gives:
\Beq\label{reduced}
\partial_t\psi\leq -\eta_{}\dt (\mu^{-1})+\mu^{-1}\,{\overline \hb}\cdot\nabla\mu-\pi\,\partial_t\rho+\csib\cdot\nabla(\partial_t\rho),
\Eeq
an inequality that replaces (\ref{dissipation}) in restricting { \emph{\`a la} Coleman \& Noll} {the possible} constitutive choices. 

On taking all of the constitutive mappings delivering $\pi,\csib,\eta$, and ${\overline \hb}$, 
depend{ent in principle} on  $\rho,\nabla\rho,\mu,\nabla\mu$, and on choosing
\Beq\label{constitutives}
\psi=\widehat\psi(\rho,\nabla\rho,\mu)=-\mu\,\rho+f(\rho)+\frac{1}{2}|\nabla\rho|^2,
\Eeq
compatibility with (\ref{reduced}) {implies that we must have:}
\Bsist \label{cn}
\begin{aligned}
\widehat\pi(\rho,\nabla\rho,\mu)&={\partial_{{\rho}} \widehat\psi(\rho,\nabla\rho,\mu)}~\displaystyle{=\mu-f'(\rho),}\\ 
\widehat\csib(\rho,\nabla\rho,\mu)&={\partial_{{\nabla\rho}} \widehat\psi(\rho,\nabla\rho,\mu)}=\nabla\rho, \\ 
\widehat\eta(\rho,\nabla\rho,\mu)&=\mu^2 \partial_{{\mu}} \widehat\psi(\rho,\nabla\rho,\mu)\displaystyle{=-\mu^2\rho,}
\end{aligned}
\Esist
together with
\Beq
\widehat{\overline \hb}(\rho,\nabla\rho,\mu,\nabla\mu) = - \widehat\Hb(\rho,\nabla\rho,\mu, \nabla\mu)\nabla \mu , \quad
\nabla \mu\cdot \widehat\Hb(\rho,\nabla\rho,\mu, \nabla\mu) \nabla\mu \geq 0 . \non
\Eeq
If {we now choose} for $\widehat\Hb$ the simplest expression $\Hb=\kappa \mathbf{1} $, implying a constant { and isotropic} mobility, and {if we once again} assume that the external distance microforce $\gamma$ and the source $\overline \sigma$ are null, {then}, with the use of (\ref{cn}) and \eqref{freeenergy}, the microforce balance (\ref{balance}) and the energy balance (\ref{energy}) become, respectively, 
\Beq\label{a}
 \div(\nabla\rho)+\mu- f'(\rho)  = 0
\Eeq
and
\Beq
 2\rho \,\dt\mu + \mu \, \dt\rho - \kappa \Delta\mu = 0,   \label{secondanew}
\Eeq
{ a nonlinear system for the unknowns $\rho$ and $\mu$} that { we supplement}
 with homogeneous Neumann condition{s at the body's boundary:}
\Beq   \label{Ibcic1}
  \dn \rho = \dn \mu = 0 
  \Eeq
(here $\dn$ denotes the outward normal derivative), 
and with the initial condition{s:}
\Beq \label{Ibcic2}
   \rho|_{t=0} = \rhoz\,,\quad\mu|_{t=0} = \mu_0\,.
\Eeq

Needless to say, \eqref{a} is the same `static' relation between $\mu $ and $\rho$ 
as $\eqref{CH}_2$. Instead, \eqref{secondanew}
is rather different from $\eqref{CH}_1${, for a number of reasons:} 
\begin{itemize}
\item \eqref{secondanew} is nonlinear 
(whereas  $\dt\rho - \kappa \Delta\mu = 0 $ is a linear equation);  
\item {the}
time derivatives of $\rho$ and $\mu$ are both present in \eqref{secondanew}; 
\item
there are {nonconstant} factors in front of both $\dt\mu $ and $\dt\rho$.
\end{itemize} 
{Moreover,  it should be possible to show that the initial/boundary-value problem \eqref{a}-\eqref{Ibcic2} has solutions $\rho\in [0,1]$ and $\mu>0$.}

{We must confess that we boldly attacked this problem as is, prompted to optimism by the successful outcome of a previous joint research effort \cite{CGPS, CGPS2}, in which we tackled mathematically the system of Allen-Cahn type one arrives at via the approach in \cite{Podio} for processes of phase segregation in the absence of diffusion. }
Unfortunately, system \accorpa{a}{Ibcic2} turned out to be  {too 
difficult for us.} {Therefore, we decided to study a regularized 
version of it,} obtained by introducing two extra terms, $\eps\, \partial_t\mu $ in 
\eqref{secondanew} and $\delta\,\partial_t\rho$ in the \lhs\ of \eqref{a}, for 
{small positive coefficients} $\eps$ and $\delta$. 

The {introduction of the} first term is motivated by the {desire to have
a strictly positive coefficient as a factor} of $\partial_t \mu$ in \eqref{secondanew},
in order to guarantee the parabolic structure of {this} equation. {As to the other term, on the one hand it gives \eqref{a}
the form of an Allen-Cahn equation with \rhs\ $\mu$; on the other hand, it assimilates our present model to the so-called \emph{viscous Cahn-Hilliard equations}} (see, e.g., 
\cite{bds, mr, Ros} and references therein). With these measures, and taking $\kappa =1 $ for simplicity, we write the following modified version of problem \accorpa{a}{Ibcic2}, with inversion of the order of the differential equations:
\Bsist
  & (\eps + 2\rho) \dt\mu + \mu \, \dt\rho - \Delta\mu = 0
  & \quad \hbox{in $\Omega \times (0,+\infty)$,}
  \label{Iprima}
  \\
  & \delta \dt\rho - \Delta\rho + f'(\rho) = \mu 
  & \quad \hbox{in $\Omega \times (0,+\infty)$,}
  \label{Iseconda}
  \\
  & \dn\mu = \dn\rho = 0
  & \quad \hbox{on $\Gamma \times (0,+\infty)$,}
  \label{Ibc}
  \\
  & \mu(\cpto,0) = \muz
  \aand
  \rho(\cpto,0) = \rhoz
  & \quad \hbox{in $\Omega$,}
  \label{Icauchy}
\Esist
\Accorpa\Ipbl Iprima Icauchy
where $\Omega \subset \erre^3$ is a bounded domain  with a sufficiently smooth boundary $\Gamma$. We {remark} that such a regularized system {has} the typical features of a phase field model, but with {a nonstandard} equation \eqref{Iprima} for the chemical potential $\mu$,
while quite often phase field systems {feature temperature} and order parameter as variables. 

By assuming, as we did in  \cite{CGPS, CGPS2}, that $f'$ is the sum of a strictly increasing $C^1$ function 
$f'_1$ with domain $(0,1)$ {that is singular at the} endpoints, and of a smooth bounded perturbation $f'_2$
(to allow for a double- or multi-well potential $f$), 
we prove the \emph{existence of a strong solution $(\mu, \rho)$ to \Ipbl\ satisfying  $\mu \geq 0$ and 
$0< \rho <1 $ almost everywhere in $\Omega \times (0,+\infty)$} (of course, the initial data have {to meet the same} requirements in $\Omega$).   
Our existence proof is rather standard; it is based on an approximation $\rightarrow$ a priori estimates $\rightarrow$ passage-to-the-limit 
procedure. Under additional assumptions, by using {certain delicate} iterative estimates, we also show that the component $\mu$ \emph{is bounded above}{; this is} probably the most difficult and technical part of the present paper. Boundedness {of $\mu$ is expedient to deduce that $f'(\rho) $ is bounded as well;} as a consequence, $\rho$ 
\emph{{stays away} from the threshold 
values $0$ and $1$}.
These boundedness properties are very useful {in proving \emph{uniqueness} of such solutions,}
since $f'(\rho)$ { can be treated as a Lipschitz-continuous} function of $\rho$. 

As a final step, we deal with the long-time behavior of the system. We prove that \emph{each element $(\mu_\omega, \rho_\omega)$ of 
the $\omega$-limit set is a steady state solution of} \Ipbl; therefore, in particular, 
$\mu_\omega$ \emph{is a constant} (cf. \eqref{Iprima} and \eqref{Ibc}). {This concludes our description of the contents of} this paper. Needless to say, {it would} be interesting and challenging 
to study the singular limit of the solutions to \Ipbl\ as $\eps$ or $\delta$ {tends to zero, or both parameters do. We plan to undertake such a study in the near future.}

%%%%%%%%%%%%%%%%%%%%%%%%%%%%%%%%%%%%%%%%%%%%%%%%%%%%%%%%%%%%%%%%%%%%%%%%

\section{Main results}
\label{MainResults}
\setcounter{equation}{0}

In this section, we 
{describe the mathematical problem under investigation}, 
make our assumptions precise, and state our results.
First of all,
we assume $\Omega$ to be a bounded connected open set in $\erre^3$
with smooth boundary~$\Gamma$
 ({to treat the lower-dimensional cases would only require} minor
changes).
Moreover, for convenience we set:
\Beq
  V := \Huno,
  \quad H := \Ldue ,
  \aand
  W := 
\graffe{v\in\Hdue:\ \dn v = 0 
\;\,{\textrm{on}\,\; \Gamma}},
  \label{defspazi}
\Eeq
and {we endow these} spaces with their standard norms,
for which we use {the self-explanato\-ry} notation 
$\normaV\cpto$ {(but $\normaH\cpto$ denotes the norm 
of any power of $H$)}.
We remark that the embeddings $W\subset V\subset H$ are compact,
because $\Omega$ is bounded and smooth.
Since $V$ is dense in~$H$,
we can identify $H$ 
{with} 
a subspace of $\Vp$ in the usual~way
(i.e.,~{ so as to have} that $_{\Vp}\<u,v>{}_V=(u,v)_H$
for every $u\in H$ and $v\in V$){; the embedding $H\subset\Vp$ is also compact.}
{As to} the potential~$f$, we assume~that
\Bsist
  \hskip-1cm && f = f_1 + f_2,
  \quad \hbox{where functions} \quad
  \hbox{$f_1,f_2:(0,1) \to \erre\;$ are { such that}}
  \label{hpf}
  \\
  \hskip-1cm && \hbox{$f_1$ is $C^1$ and convex}, \quad
  \hbox{$f_2$ is $C^2$},
  \quad
  \hbox{$f_2''$ is bounded},
  \label{hpfi}
  \\
  \hskip-1cm && \lim_{r\searrow0} f_1'(r) = - \infty\,,
  \aand
  \lim_{r\nearrow1} f_1'(r) = + \infty .
  \label{hpfp}
\Esist
\Accorpa\Hpf hpf hpfp
For the initial data, we stipulate that
\Bsist
  && \muz \in V
  \aand \muz \geq 0 \quad \aeO;
  \label{hpmuz}
  \\
  && \rhoz \in W, \quad
  0 < \rhoz < 1 \quad \hbox{in $\Omega$};
  \aand
  f'(\rhoz) \in H .
  \label{hprhoz}
\Esist
\Accorpa\Hpdati hpmuz hprhoz
We stress that 
{the conditions in}
\eqref{hprhoz} imply
{that}
\Beq
  \rhoz \in \Cx0
  \aand
  f(\rhoz) \in H.
  \label{regrhoz}
\Eeq
Indeed,
$W\subset\Cx0$,
{assumptions}
\eqref{hpfi} hold, and,
{by convexity,}
$-c\leq f_1(\rhoz)\leq f_1(1/2)+f_1'(\rhoz)(\rhoz-1/2)$
for some $c\in\erre$.

Our aim is 
{to solve}
problem~\Ipbl\ in a strong sense,
i.e., we want to find a pair $(\mu,\rho)$
of { such smooth} functions { satisfying suitable summability conditions and
unilateral constraints that \Ipbl\ are made fully} meaningful.
{Precisely}, we fix a final time
{$T>0$,  we set $Q:=\Omega\times(0,T)$, and we} require~that:
\Bsist
  && \mu \in \H1H \cap \L2W,
  \label{regmu}
  \\
  && \rho \in \W{1,\infty}H \cap \H1V \cap \L\infty W,
  \label{regrho}
  \\
  && \mu \geq 0 \quad \aeQ,
  \label{mupos}
  \\
  && 0 < \rho <1 \quad \aeQ
  \aand
  f'(\rho) \in \L\infty H.
  \label{regfprho}
\Esist
\Accorpa\Regsoluz regmu regfprho
%
%{This}
%notation is widely used in the 
%{following}.
Note that the boundary conditions~\eqref{Ibc} follow from
\accorpa{regmu}{regrho}, due to the definition of~$W$
in~\eqref{defspazi}.
In conclusion, we look for $(\mu,\rho)$ satisfying~\Regsoluz\
and fulfilling the system
\Bsist
  & (\eps + 2\rho) \dt\mu + \mu \, \dt\rho - \Delta\mu = 0
  & \quad \aeQ,
  \label{prima}
  \\
  & \delta \dt\rho - \Delta\rho + f'(\rho) = \mu
  & \quad \aeQ,
  \label{seconda}
  \\
  & \mu(0) = \muz
  \aand
  \rho(0) = \rhoz
  & \quad \aeO .
  \label{cauchy}
\Esist
\Accorpa\Pbl prima cauchy
Here is our {main} result.

\Bthm
\label{Esistenza}
Assume 
{that}
\Hpf\ and \Hpdati
{\,are satisfied}.
Then, there exists a pair $(\mu,\rho)$ satisfing \Regsoluz\ and solving
problem~\Pbl.
\Ethm

{Once existence is secured}, one wonders about uniqueness.
{We are able to prove it for solutions having the following additional properties:}
\Beq
  \mu \in \LQ\infty ; \quad
  \inf\rho >0 \aand \sup\rho <1. 
  \label{regsoluzbis}
\Eeq

\Bthm
\label{Unicita}
Assume
{that}
\Hpf\ and \Hpdati
{\,are satisfied}.
Then, any two solutions to problem \Pbl\
satisfing \Regsoluz\ and~\eqref{regsoluzbis}
coincide.
\Ethm
\noindent{Interestingly, t}he additional boundedness 
{conditions for
$(\mu,\rho)$ postulated above
are fulfilled}
whenever the data of the problem 
{have similar boundedness properties}, in addition to \Hpdati.
%{Indeed, we have the following result.}

\Bthm
\label{Limitatezza}
Assume 
{that}
{\accorpa{hpf}{hprhoz}} and
{the following conditions are satisfied:} 
\Beq
  \muz \in \Lx\infty; \quad
  \inf\rhoz >0 \aand \sup\rhoz <1 .
  \label{hpdatibis}
\Eeq
Then, any pair $(\mu,\rho)$ satisfing \Regsoluz\ and solving
problem~\Pbl\ 
satisfies \eqref{regsoluzbis} {as well}.
\Ethm

\Brem
\label{Piuregolare}
Even though the regularity of the solution
given by \Regsoluz\ and~\eqref{regsoluzbis}
is completely satisfactory for our purposes,
we observe that some further smoothness can be proved
once the properties \eqref{regsoluzbis} are established.
In that case, equation~\eqref{seconda} can be read in the form
$\dt\rho-\Delta\rho=g$ with $g\in\LQ\infty$,
whence further regularity for $\rho$ can be derived,
and a bootstrap procedure can start.
Indeed, further regularity for $\rho$ implies
that stronger properties for $\mu$ can be proved by~\eqref{prima}.
This improves the regularity of $g$ and leads to 
{an increase of}
the regularity of~$\rho$.
\Erem

Once \wepo ness on every finite time interval is ensured, 
one can study the \loti\ \bhv\ of the solution.
In particular, one can try to \characteriz e the $\omega$-limit
of any trajectory $(\mu,\rho)$ in some topology.
We choose the weak topology of~$H\times V$
and define such an $\omega$-limit as follows:
\Bsist
  & \omega(\mu,\rho)
  & = \bigl\{
    (\muo,\rhoo):\ \bigl( \mu(t_n),\rho(t_n) \bigr) \to (\muo,\rhoo)\ 
  \non
  \\
  && \qquad
    \hbox{weakly in $H\times V$ for some sequence $t_n\mathrel{\scriptstyle\nearrow}+\infty$}
  \bigr\} .
  \label{defomegalim}
\Esist
Our last result gives a relationship between such an $\omega$-limit
and the set of \emph{steady states}, i.e.,
the set of the time-independent solutions $(\mus,\rhos)$ to
\accorpa{prima}{seconda} with homogeneous Neumann boundary condition
satisfying natural regularity properties.
Note that in such a case $\mus$~must be harmonic, thus constant,
since $\Omega$ is connected.
Therefore, a~steady state is a pair $(\mus,\rhos)$
{such that}
$\mus$~is a 
{nonnegative}
constant and $ \rhos$ solves the following problem:
\Beq
  \rhos \in W , \quad
  0 < \rhos < 1 , \quad
  f'(\rhos) \in H ,
  \aand
  -\Delta\rhos + f'(\rhos) = \mus
  \quad \aeO 
  \label{secondas}
\Eeq
(there is no reason for $\rhos$ to be constant,
since $f$ is not required to be convex).
%Here is the result we have.

\Bthm
\label{Longtime}
Assume 
{that conditions}
\Hpf, \Hpdati, and \eqref{hpdatibis},
{are satisfied}.
Let $(\mu,\rho)$ be the corresponding solution
satisfing \Regsoluz\ and~\eqref{regsoluzbis}.
Then, the $\omega$-limit $\omega(\mu,\rho)$
is
{nonempty}, 
compact, and connected in the weak topology of $H\times V$; moreover,
each of its elements coincides with a steady state $(\mus,\rhos)$ (that is to say,
$\mus$~is a 
{nonnegative}
constant
and $\rhos$~solves~\eqref{secondas}).
\Ethm

Our paper is \organiz ed as follows.
In the next section, we prove Theorem~\ref{Esistenza},
while Theorems~\ref{Unicita} and~\ref{Limitatezza}
are
proved in Section~\futuro\ref{UniqBdd}.
Our last section is devoted to the proof of Theorem~\ref{Longtime}.

%\medskip
%\Brem 

Throughout the paper,
%\marginpar{\footnotesize{ "throughout the whole" 
%is a bit redundant!}}
we account for the \wk\ embedding $V\subset\Lx q$
for $1\leq q\leq 6$
and the related Sobolev inequality:
\Beq
  \norma v_{\Lx q} \leq C \normaV v
  \quad \hbox{for every $v\in V$ and $1\leq q \leq 6$,}
  \label{sobolev}
\Eeq
where $C$ depends on~$\Omega$ only,
since sharpness is not needed (the embedding $V\subset\Lx q$ is
compact if $q<6$). 
Furthermore, we {repeatedly make use of}
the \wk\ \holder\ inequality,
the interpolation inequality
\Bsist
  \hskip-1cm && \norma v_{\Lx r}
  \leq \norma v_{\Lx p}^\theta \, \norma v_{\Lx q}^{1-\theta}
  \quad \hbox{for $v\in\Lx p\cap\Lx q$},
  \non
  \\
  \hskip-1cm && \qquad \hbox{where} \quad
  p,q,r \in [1,+\infty] ,
  \quad \theta \in [0,1] ,
  \aand
  \frac 1r = \frac \theta p + \frac {1-\theta} q,
  \label{Interpolazione}
\Esist
%\Accorpa\Interpolazione interpol paraminterpol
and the elementary Young inequality
\Beq
  ab \leq \sigma a^2 + \frac 1{4\sigma} \, b^2
  \quad \hbox{for every $a,b\geq 0$ and $\sigma>0$}.
  \label{young}
\Eeq
Finally, throughout the paper
%We conclude the present section by stating a general rule
%that we use as far as constants are concerned,
%in order to avoid a boring notation.
we use a {small-case italic} $c$ for different constants, that
{may only} depend 
on~$\Omega$, the final time~$T$, the shape of~$f$, 
the properties of the data involved in the statements at hand,
and the coefficients $\eps$ and~$\delta$; {a~notation like~$c_\sigma$ 
signals a constant that depends also on the parameter~$\sigma$}. {The reader should keep in mind that} the meaning of $c$ and $c_\sigma$ might
change from line to line and even in the same chain of inequalities, 
{whereas those constants we need to refer to are always denoted by 
capital letters, just like $C$ in~\eqref{sobolev}.}
%\Erem
%%%%%%%%%%%%%%%%%%%%%%%%%%%%%%%%%%%%%%%%%%%%%%%%%%%%%%%%%%%%%%%%%%%%%%%%

\section{Existence}
\label{Existence}
\setcounter{equation}{0}

In this section, we prove Theorem~\ref{Esistenza}.
Our method uses an 
{approximation} scheme based on a time delay in the \rhs\
of~\eqref{seconda}.
Namely, we define the translation operator $\T:\L1H\to\L1H$ depending
on a time step $\tau>0$ by setting, for $v\in\L1H$ and \aat,
\Beq
  (\T v)(t) := v(t-\tau)
  \quad \hbox{if $t>\tau$}
  \aand
  (\T v)(t) := \muz
  \quad \hbox{if $t<\tau$,}
  \label{defT}
\Eeq
and consider the problem obtained by replacing the \rhs\
of~\eqref{seconda} by~$\T\mu$, i.e.,
we look for a pair $(\mut,\rhot)$ such~that
\Bsist
  & \hbox{$(\mut,\rhot)$ satisfies \Regsoluz}
  \label{regsoluztau}
  \\
  & (\eps + 2 \rhot) \dt\mut - \Delta\mut + \mut \, \dt\rhot = 0
  & \quad \aeQ
  \label{primatau}
  \\
  & \delta \dt\rhot - \Delta\rhot + f'(\rhot) = \T\mut
  & \quad \aeQ
  \label{secondatau}
  \\
  & \mut(0) = \muz
  \aand
  \rhot(0) = \rhoz
  & \quad \aeO .
  \label{cauchytau}
\Esist
\Accorpa\Pbltau regsoluztau cauchytau
For convenience, we allow $\tau$ to take just discrete values,
namely, $\tau=T/N$, where $N$ is any positive integer.
Our existence proof consists 
{in}
two parts.
Firstly, we check 
that problem~\Pbltau\ is \wepo\,(see {the next lemma). Secondly,}
we let $\tau$ tend to~$0$.
This is done by proving a number of a~priori
estimates and using compactness and monotonicity { arguments}.

\Blem
\label{Existtau}
There exists a unique pair $(\mut,\rhot)$ solving problem \Pbltau.
\Elem

\par\noindent{\bf Proof.}
{Recall} that $\tau=T/N$.
Hence, if we set $t_n:=n\tau$ for $n=0,\dots,N$,
we see that problem \Pbltau\ {becomes} equivalent to a finite sequence of $N$ problems
that can be solved step by step.
However, instead of considering the natural time intervals
$[t_{n-1},t_n]$, $n=1,\dots,N$, and glueing the solutions 
{together},
we solve $N$ problems
on the time intervals $I_n=[0,t_n]$, $n=1,\dots,N$,
by constructing the solution directly on the whole of~$I_n$ at each step.
{These} problems are the following:
\Bsist
  && (\eps + 2 \rhon) \dt\mun - \Delta\mun  + (\dt\rhon) \mun = 0
  \aand
  \mun \geq 0
  \quad \hbox{a.e.\ in $\Omega\times I_n$}
  \label{priman}
  \\
  && \dn\mun(t)|_\Gamma = 0 \quad \hbox{for a.a.\ $t\in I_n$}
  \aand
  \mun(0) = \muz
  \label{ibcmun}
  \\
  && 0 < \rhon < 1
  \aand
  \delta \dt\rhon - \Delta\rhon + f'(\rhon) = \T\munmu
  \quad \hbox{a.e.\ in $\Omega\times I_n$}
  \label{secondan}
  \\
  && \dn\rhon(t)|_\Gamma = 0 \quad \hbox{for a.a.\ $t\in I_n$}
  \aand
  \rhon(0) = \rhoz.
  \label{ibcrhon}
\Esist
{Their}
 solutions are required to satisfy the regularity properties induction
obtained by 
{taking}
$t_n$ in place of $T$ in~\Regsoluz.
The operator $\T$ that appears on the \rhs\ of~\eqref{secondan} 
acts on functions that are not defined in the whole of~$(0,T)$.
However, its meaning is still given by~\eqref{defT} if $n>1$,
while we simply set $\T\munmu=\muz$ if $n=1$.

Clearly, the solution $(\mut,\rhot)$ we are looking for
is simply given by $(\mu_N,\rho_N)$.
The above problems can be solved inductively, {because}
the \rhs\ of \eqref{secondan} is known at each step {{in the next lemma}:}
{one first solves}
problem 
\accorpa{secondan}{ibcrhon} for~$\rhon$,
and then problem \accorpa{priman}{ibcmun} for~$\mun$.
We note that the former {problem} is quite standard; 
the latter is a regular linear parabolic problem (the~coefficient of $\dt\mun$ is~$\geq\eps$)
provided that $\dt\rhon$ is 
{sufficiently}
smooth. {That} the inequality $\mun\geq0$ {holds} is not obvious.
The uniqueness of a solution $(\mun,\rhon)$
satisfying smoothness properties analogous to~\Regsoluz\ is clear,
and the existence of a variational solution is expected.
However, 
{not even the desired regularity is obviously guaranteed. Therefore, we provide a few arguments in this direction}. 

{Proceeding at an as-low-as-possible level of formality, we introduce 
a problem depending on a positive parameter~$\lambda$ and approximating problem \accorpa{secondan}{ibcrhon}. To begin with, we}
%{To this end, we}
 \regulariz e $f_1$ and $f_2$
by constructing certain suitable $C^2$~approximations  $\ful,\fdl$
{having bounded first and second derivatives. Precisely,} we assume that $\fdl''$ is bounded uniformly with respect to~$\lambda$; moreover, on thinking of $f_1'$ as a maximal monotone graph in~$\erre\times\erre$, 
we assume that $\ful$ is convex and that $\ful'$
is similar to the Yosida \regulariz ation of~$f_1'$
(see, e.g., \cite[p.~28]{Brezis}),
in order to preserve the main properties of the latter
(such a \regulariz ation is detailed, e.g., in~\cite[Section~3]{GR}). Finally, we set $\fl=\ful+\fdl$.
The approximating problem~is: 
%{given by}
\Bsist
  && \delta \dt\rhonl - \Delta\rhonl + \fl'(\rhonl) = \T\munmu
  \quad \hbox{a.e.\ in $\Omega\times I_n$,}
  \label{secondanl}
  \\  
  && \dn\rhonl(t)|_\Gamma = 0 \quad \hbox{for a.a.\ $t\in I_n$},
  \aand
  \rhonl(0) = \rhoz;
  \label{ibcrhonl}
\Esist
it has a unique smooth solution,
which satisfies 
{sufficiently strong}
a~priori estimates {to allow 
 letting} $\lambda$ tend to zero in \accorpa{secondanl}{ibcrhonl}.
This leads to a solution $\rhon$ to problem~\accorpa{secondan}{ibcrhon},
{which can be used} to solve problem~\accorpa{priman}{ibcmun}.
{Needless to say,} the desired regularity for $\rhon$
{will follow once we prove}
suitable estimates uniformly with respect to~$\lambda$.
We confine ourselves to derive the highest-order estimate{,
the others being} quite standard.

{With a view toward assembling a proof by} induction, we assume~that
\Beq
  \munmu \in H^1(I_{n-1};H) \cap L^\infty(I_{n-1};V)
  \aand
  \munmu \geq 0 \quad{\textrm{for}\quad n>1}
  \label{hpinduz}
\Eeq
%{(not requiring anything else in the case $n=1$)}
%, which makes $I_{n-1}$ thus \eqref{hpinduz} meaningless)},
% \marginpar{\footnotesize{mi sembra detto male e, comunque, non mi pare che 
% si sia detto che $I_1=[0,t_1]=\{0\}$. Sarei propenso a sopprimere tutta la 
% parentesi, anche in vista di quello che vien detto subito dopo.}}
and we prove that
\Bsist
  && \norma{\rhonl}_{W^{1,\infty}(I_n;H) \cap H^1(I_n;V) \cap L^\infty(I_n;W)}
  + \norma{\ful'(\rhonl)}_{L^\infty(I_n;H)}
  \leq c_\tau,
  \label{stimarhon}
  \\
  && \mun \in H^1(I_n;H) \cap L^\infty(I_n;V) \cap L^2(I_n;W)
  \aand
  \mun \geq 0 
  \label{stimamun}
\Esist
{(as anticipated in closing Section 2, in \eqref{stimarhon} as well as in the 
following 
the symbol $c_\tau$ stands for one or another of a list of different constants
that do not depend on~$\lambda$, but are allowed to depend on~$\tau$).}
We remark that the induction procedure can actually start,
{because} $\T\munmu=\muz$ if $n=1$ and, moreover,
%{the}
properties~\eqref{hpmuz} and~\eqref{hprhoz} for $\muz$ and $\rhoz$ 
{are fulfilled; these properties are also used at each step.}

{We omit  stressing the dependences on $n$ and $\lambda$, and write simply $u$ and $\uz$ for, respectively, $\dt\rhonl$ and $\dt\rhonl(0)$}.
By differentiating~\eqref{secondanl} with respect to time,
we see that $u$ solves the equation:
\Beq
  \delta \dt u - \Delta u + \ful''(\rhonl) \, u
  = \dt(\T\munmu) - \fdl''(\rhonl) \, u
  \quad \hbox{a.e.\ in $\Omega\times I_n$},
  \label{dtsecondanl}
\Eeq
and satisfies {both} the Cauchy condition
$u(0)=\uz$
and homogeneous Neumann boundary condition.
Hence, by testing~\eqref{dtsecondanl} by~$u$
and using the convexity of~$\ful$, we immediately obtain
for $t\in I_n$ 
{that}
\Beq
  \frac\delta 2 \, \normaH{u(t)}^2
  + \intQt |\nabla u|^2
  \leq \frac\delta 2 \, \normaH{\uz}^2
  + \bigl( 1 + \sup |\fdl''| \bigr) \intQt u^2
  + \norma{\dt(\T\munmu)}_{L^2(I_n;H)}^2 .
  \label{dadtsecondanl}
\Eeq
Now, we observe that the last norm is finite, {in view of} our assumption~\eqref{hpinduz},
and that $|\fdl''|\leq c$. 
Moreover,
{due to~\eqref{secondanl}, we have that~} 
$\delta \uz=\muz+\Delta\rhoz-\fl'(\rhoz)$.
Hence, $\uz$~is bounded in~$H$, by~\accorpa{hpmuz}{hprhoz} 
and our choice of the approximation $\fl$ of~$f$.
Therefore, thanks to the Gronwall lemma, we obtain:
\Beq
  \norma u_{L^\infty(I_n;H) \cap L^2(I_n;V)} \leq c_\tau \,,
  \quad \hbox{whence} \quad
  \norma\rhonl_{W^{1,\infty}(I_n;H) \cap H^1(I_n;V)} \leq c_\tau \,.
  \non
\Eeq
Next, coming back to~\eqref{secondanl},
we deduce that $-\Delta\rhonl+\ful'(\rhonl)$ is bounded in $L^\infty(I_n;H)$ and {hence, by
a standard argument ({for instance, by testing}~\eqref{secondanl} by~$\ful'(\rhonl)$), that each of  $-\Delta\rhonl$ and $\ful'(\rhonl)$ is bounded. With this, given that} 
{the $W$-estimate follows from elliptic theory,
\eqref{stimarhon}~is established,} and we can let $\lambda$ tend to zero.
We obtain:
\Beq
  \rhon \in W^{1,\infty}(I_n;H) \cap H^1(I_n;V) \cap L^\infty(I_n;W) , \quad
  0 < \rhon < 1 ,
  \aand
  f_1'(\rhon) \in L^\infty(I_n;H) .
  \non
\Eeq
At this point, we should prove~\eqref{stimamun}.
However, we confine ourselves to derive a formal estimate
that clearly shows that the desired regularity for $\mun$
can be deduced by \regulariz ing the linear problem \accorpa{priman}{ibcmun}
(if~the coefficient $\dt\rhon$ is replaced by a smooth function and the initial datum is \regulariz ed,
{by the same token} $\dt\mun$ is an admissible test function).
For convenience, we write \eqref{priman} in the form:
\Beq
  (\eps+2\rhon) \dt\mun + \mun - \Delta\mun
  = (1-\dt\rhon) \mun;
  \non
\Eeq
next, we multiply this relation by $\dt\mun$ and use the 
{result}
 in the calculation given below.
{Since $\rhon\geq0$, we find, for $t\in I_n$,}
\Bsist
  && \eps\intQt |\dt\mun|^2
  + \frac 12 \, \normaV{\mun(t)}^2
  \leq \intQt (\eps + 2 \rhon) |\dt\mun|^2
  + \frac 12 \, \normaV{\mun(t)}^2
  \non
  \\
  && = \intQt (\eps + 2 \rhon) |\dt\mun|^2
  + \frac 12 \, \normaV\muz^2
  + \frac 12 \intQt \dt \bigl( |\mun|^2 + |\nabla\mun|^2 \bigr)
  \non
  \\
  && = \frac 12 \, \normaV\muz^2
  + \intQt \bigl(
    (\eps + 2 \rhon) |\dt\mun|^2
    + \mun \, \dt\mun + \nabla\mun \cdot \nabla\dt\mun
  \bigr)
  \non
  \\
  && = \frac 12 \, \normaV\muz^2
  + \intQt (1-\dt\rhon) \, \mun \, \dt\mun
  \non
  \\
  && \leq \frac 12 \, \normaV\muz^2
  + \iot \norma{1+|\dt\rhon(s)|}_{\Lq} \,
         \norma{\mun(s)}_{\Lq} \,
         \norma{\dt\mun
         {(s)}
         }_{\Lx2} \, ds
  \non
  \\
  && \leq \frac 12 \, \normaV\muz^2
  + \frac \eps 2 \intQt |\dt\mun|^2
  + \frac {C^2}{2\eps}
    \iot \normaV{1+|\dt\rhon(s)|}^2 \, \normaV{\mun(s)}^2 \, ds,
  \non
\Esist
by the Sobolev and Young inequalities~\eqref{sobolev} and~\eqref{young}.
Then, the Gronwall lemma yields that
\Beq
  \norma{\dt\mun}_{L^2(I_n;H)} + \norma\mun_{L^\infty(I_n;V)} \leq c_M, 
  \label{prestimamun}
\Eeq
where $M$ is a constant satisfying
$M\geq\normaV\uz+\norma{\dt\rhon}_{L^2(I_n;V)}$.
By comparison in~\eqref{secondan},
even $\Delta\mun$ is estimated in~$L^2(I_n;H)$,
since a bound for $\mun\dt\rhon$ in the same space follows from~\eqref{prestimamun}.
By elliptic regularity, we derive the desired estimate for $\mun$ in~$L^2(I_n;W)$.
So, the first 
{assertion in} \eqref{stimamun} is established; it remains for us to show that $\mun\geq0$.
This is done by testing \eqref{priman} by~$-\mun^-$.
We obtain, for $t\in I_n$, that
\Bsist
  && \frac 12 \intQt \dt \bigl( (\eps + 2 \rhon) |\mun^-|^2 \bigr)
  + \intQt |\nabla \mun^-|^2
  \non
  \\
  && = \intQt \bigl(
    (\eps+2\rhon) \dt\mun (-\mun^-)
    + (\dt\rhon) \mun (-\mun^-)
    + \nabla\mun \cdot \nabla(-\mun)^-
  \bigr)
  = 0 .
  \non
\Esist
As $\rhon\geq0$ and $\muz\geq0$, we deduce that
\Beq
  \eps \iO |\mun^-(t)|^2
  \leq \iO (\eps + 2 \rhon(t)) |\mun^-(t)|^2 
  \leq \iO (\eps + 2 \rhoz) |\muz^-|^2
  = 0,
  \non
\Eeq
whence 
{it immediately follows that}
$\mun^-=0$, i.e., that $\mun\geq0$.
Thus, the lemma is proved.
\QED

{Now that}
the \wepo ness of problem \Pbltau\ is established,
we perform a number of a~priori estimates of its solution. {These estimates allow us to let $\tau$ tend to zero, so as to}
%{thus}
prove our existence result for problem~\Pbl.
In order to make the formulas {to come more readable, 
we shall omit the index~$\tau$ in the calculations, waiting for writing} $(\mut,\rhot)$ only when each estimate is established.

\step First a priori estimate

We observe that
$\dt\bigl((\eps/2)\mu^2+\rho\mu^2\bigr)=\bigr((\eps+2\rho)\dt\mu+\mu\,\dt\rho\bigr)\mu$.
Thus, 
{testing \eqref{primatau} by $\mu$ and integrating, we} 
obtain, for $t\in(0,T)$, that
\Beq
  \iO \Bigl( \frac \eps 2 \, \mu^2 + \rho \mu^2 \Bigr)(t)
  + \intQt |\nabla\mu|^2
  = \iO \Bigl( \frac \eps 2 \, \muz^2 + \rhoz \muz^2 \Bigr)
  = c .
  \label{perprimastima}
\Eeq
This implies that
\Beq
  \norma\mut_{\L\infty H\cap\L2V} \leq c .
  \label{primastima}
\Eeq

\step Second a priori estimate

This standard estimate for phase field equations
can be derived by testing \eqref{secondatau} by~$\dt\rho$.
We~get:
\Beq
  \norma\rhot_{\H1H\cap\L\infty V}
  + \norma{f(\rhot)}_{\L\infty\Luno}
  \leq c .
  \label{secondastima}
\Eeq

\step Third a priori estimate

We rewrite \eqref{secondatau} as
\Beq
  -\Delta\rho + f_1'(\rho) = - \delta\dt\rho - f_2'(\rho) + \T\mu,
  \label{venti}
\Eeq
and notice that the \rhs\ is bounded in~$\L2H$.
Then, by applying a standard procedure
({for instance, testing}
 by $f_1'(\rho)$),
and {counting on} elliptic regularity, we deduce~that
\Beq
  \norma\rhot_{\L2W} + \norma{f_1'(\rhot)}_{\L2H} \leq c .
  \label{terzastima}
\Eeq

\step Fourth a priori estimate

{To derive the next inequality, we prefer to proceed  formally, 
avoiding the $\lambda$-\regulariz ation we used in the proof of Lemma~\ref{Existtau}.}
As  {for \eqref{dadtsecondanl}}, we~obtain
{the following estimate:}
\Bsist
  && \frac\delta 2 \, \normaH{\dt\rho(t)}^2
  + \intQt |\nabla\dt\rho|^2
  \non
  \\
  && \leq \frac\delta 2 \, \normaH{\Delta\rhoz-f_1'(\rhoz)+\muz}^2
  + \sup |\fdl''| \intQt |\dt\rho|^2
  + \intQt (\dt\T\mu) \, \dt\rho .
  \label{perquartastima}
\Esist
Once 
{this}
inequality is established,
{our} procedure is rigorous.
The estimate of the last term requires {now} more care than before,
{because we aim to} obtain bounds that are uniform with respect to~$\tau$.
We~have:
\Beq
  \intQt (\dt\T\mu) \, \dt\rho 
  = \int_\tau^t \iO  \dt\mu(s-\tau) \, \dt\rho(s) \, ds
  = \iotmt \iO  \dt\mu(s) \, \dt\rho(s+\tau) \, ds,
  \non
\Eeq
and we compute $\dt\mu$ from~\eqref{priman}.
On recalling that $\rho\geq0$, we can continue as follows:
\Bsist
  && \intQt (\dt\T\mu) \, \dt\rho 
  = \iotmt \iO \Bigl(
      \frac 1 {\eps+2\rho} \, \bigl( \Delta\mu - \mu \, \dt\rho \bigr)
    \Bigr)(s) \, \dt\rho(s+\tau)\, {ds}
  \non
  \\
  && = \iotmt \iO \Bigl\{
       - \Bigl( \frac {\nabla\mu} {\eps+2\rho} \Bigr)(s) \cdot \nabla\dt\rho(s+\tau)
       + 2 \, \frac {\dt\rho(s+\tau)} {(\eps+2\rho(s))^2} \, \nabla\mu(s) \cdot \nabla\rho(s)
  \non
  \\
  && \qquad \qquad \qquad {}- \dt\rho(s) \, \mu(s) \, \dt\rho(s+\tau) \, \frac 1 {\eps+2\rho(s)} \,
     \Bigr\} \,{ds}
  \non
  \\
  && \leq \frac 14 \intQt |\nabla\dt\rho|^2
  + c \norma\mu_{\L2V}^2
  \non
  \\
  && \quad {} + c \iotmt \norma{\dt\rho(s+\tau)}_{\Lq} \, \normaH{\nabla\mu(s)} \, \norma{\nabla\rho(s)}_{\Lq} \, ds
  \non
  \\
  && \quad {} + c \iotmt \norma{\dt\rho(s+\tau)}_{\Lq} \, \norma{\mu(s)}_{\Lq} \, \normaH{\dt\rho(s)} \, ds. 
  \label{dtTmu}
\Esist
{We need to}
estimate the last two integrals.
{As to the first,
we begin by using} the Sobolev inequality~\eqref{sobolev}
and the elementary Young inequality~\eqref{young}.
We~find that
\Bsist
  && \iotmt \norma{\dt\rho(s+\tau)}_{\Lq} \, \normaH{\nabla\mu(s)} \, \norma{\nabla\rho(s)}_{\Lq} \, ds
  \non
  \\
  && \leq \frac 18 \iotmt \normaV{\dt\rho(s+\tau)}^2 \, ds
  + c \iot \normaH{\nabla\mu(s)}^2 \, \normaV{\nabla\rho(s)}^2 \, ds
  \non
  \\
  && \leq \frac 18 \intQt |\nabla\dt\rho|^2 
  + \frac 18 \, \norma{\dt\rho}_{\L2H}^2 
  \non
  \\
  && \quad {}
  + c \iot \normaV{\mu(s)}^2 \, \bigl( \normaV{\rho(s)}^2 + \normaH{\Delta\rho(s)}^2 \bigr) \, ds,
  \non
\Esist
the last inequality {holding because, thanks to elliptic regularity,}
$\normaW v\leq c(\normaV v+\normaH{\Delta v})$
for 
{any}
$v\in V$ such that $\Delta v\in H$ and $\dn v|_\Gamma=0$.
At this point, we recall that $\rho$ is bounded in $\H1H\cap\L\infty V$, and
%{that}
 $\mu$  in $\L2V$,
by \eqref{secondastima} and~\eqref{primastima}.
{Moreover,}
{we}
notice that \eqref{venti} entails (formally, by testing it by~$-\Delta\rho(s)$):
\Bsist
  \normaH{\Delta\rho(s)}^2
  &\leq & \delta^2 \normaH{\dt\rho(s)}^2
  + c \bigl( 1 + \normaH{\T\mu(s)}^2 \bigr) \non \\[0.2cm]
  &\leq & \delta^2 \normaH{\dt\rho(s)}^2  + c
  \quad \hbox{for a.a.\ $s\in(0,T)$},
  \non
\Esist
the last inequality being a consequence of~\eqref{primastima}.
Therefore, we 
{can} 
infer that
\Bsist
  && \iotmt \norma{\dt\rho(s+\tau)}_{\Lq} \, \normaH{\nabla\mu(s)} \, \norma{\nabla\rho(s)}_{\Lq} \, ds
  \non
  \\[0.2cm]
  && \leq \frac 18 \intQt |\nabla\dt\rho|^2
  + c + c \iot \normaV{\mu(s)}^2 \, \normaH{\dt\rho(s)}^2 \, ds.
  \label{penultimo}
\Esist
{Passing now to estimate the} last integral in~\eqref{dtTmu},
we have that
\Bsist
  && \iotmt \norma{\dt\rho(s+\tau)}_{\Lq} \, \norma{\mu(s)}_{\Lq} \, \normaH{\dt\rho(s)} \, ds 
  \non
  \\
  && \leq \frac 18 \iotmt \normaV{\dt\rho(s+\tau)}^2 \, ds
  + c \iot \normaV{\mu(s)}^2 \, \normaH{\dt\rho(s)}^2 \, ds
  \non
  \\
  && \leq \frac 18 \intQt |\nabla\dt\rho|^2
  + \frac 18 \, \norma{\dt\rho}_{\L2H}^2
  + c \iot \normaV{\mu(s)}^2 \, \normaH{\dt\rho(s)}^2 \, ds
  \non
  \\
  && \leq \frac 18 \intQt |\nabla\dt\rho|^2
  + c
  + c \iot \normaV{\mu(s)}^2 \, \normaH{\dt\rho(s)}^2 \, ds .
  \non
\Esist
With this and~\eqref{penultimo},
we see that \eqref{dtTmu} yields:
\Beq
  \intQt (\dt\T\mu) \, \dt\rho 
  \leq c + \frac 12 \intQt |\nabla\dt\rho|^2
  + c \iot \normaV{\mu(s)}^2 \, \normaH{\dt\rho(s)}^2 \, ds,
  \non
\Eeq
so that \eqref{perquartastima} takes the form:
\Beq
  {\frac\delta 2 \, \normaH{\dt\rho(t)}^2
  + \frac 12 \intQt |\nabla\dt\rho|^2
  \leq c  + c \iot \normaV{\mu(s)}^2 \, \normaH{\dt\rho(s)}^2 \, ds.}
  \non
\Eeq
{Since $\mu$ has been estimated in~$\L2V$, the Gronwall
lemma can be applied, so as to have that} 
\Beq
  \norma{\dt\rhot}_{\L\infty H\cap\L2V} \leq c. 
  \label{quartastima}
\Eeq
{Finally, the same argument as in the derivation of~\eqref{terzastima} yields that}
\Beq
  \norma\rhot_{\L\infty W} + \norma{f_1'(\rhot)}_{\L\infty H} \leq c .
  \label{daquartastima}
\Eeq

\step Fifth a priori estimate

We formally test \eqref{primatau} by $\dt\mu$, and obtain:
\Bsist
  && \eps \intQt |\dt\mu|^2
  + \frac 12 \iO |\nabla\mu(t)|^2
  \non
  \\  
  && \leq \frac 12 \, \normaH{\nabla\muz}^2
  - \intQt \dt\rho \, \mu \, \dt\mu
  \non
  \\
  && \leq c
  + \frac 1 {2\eps} \iot \norma{\dt\rho(s)}_{\Lq}^2 \, \norma{\mu(s)}_{\Lq}^2 \, ds
  + \frac \eps 2 \iot \normaH{\dt\mu(s)}^2 \, ds
  \non
  \\
  && \leq c
  + c \iot \normaV{\dt\rhos(s)}^2 \, \bigl( \normaH{\mu(s)}^2 + \normaH{\nabla\mu(s)}^2 \bigr) \, ds
  + \frac \eps 2 \iot \normaH{\dt\mu(s)}^2 \, ds
  \non
  \\
  && \leq c
  + c \iot \normaV{\dt\rhos(s)}^2 \, \normaH{\nabla\mu(s)}^2 \, ds
  + \frac \eps 2 \iot \normaH{\dt\mu(s)}^2 \, ds, 
  \non
\Esist
{where the last inequality follows from \eqref{quartastima} and~\eqref{primastima}.
Using}
\eqref{quartastima} once more, we can apply the Gronwall lemma
and conclude~that
\Beq
  \norma\mut_{\H1H\cap\L\infty V} \leq c .
  \label{quintastima}
\Eeq

\step Sixth a priori estimate

Recalling that $0<\rho<1$, and using the Sobolev inequality~\eqref{sobolev}, we~get:
\Bsist
  && \norma{(\eps+2\rho)\dt\mu+\mu\,\dt\rho}_{\L2H}
  \non
  \\
  && \leq (\eps+2) \norma{\dt\mu}_{\L2H} + \norma\mu_{\L\infty{\Lq}} \, \norma{\dt\rho}_{\Ldlq}
  \non
  \\
  && \leq c \bigl( \norma{\dt\mu}_{\L2H} + \norma\mu_{\L\infty V} \, \norma{\dt\rho}_{\L2V} \bigr).
  \non
\Esist
{Since}
the \rhs\ is bounded by \eqref{quintastima} and~\eqref{quartastima},
a~comparison in \eqref{primatau} shows that $\Delta\mu$ is bounded in~$\L2H$.
Consequently,  by elliptic regularity we deduce that
\Beq
  \norma\mut_{\L2W} \leq c.
  \label{sestastima}
\Eeq

\step Conclusion

Collecting all the estimates we have proved, 
we see~that
\Bsist
  && \mut \to \mu
  \quad \hbox{weakly star in $\H1H\cap\L\infty V\cap\L2W$},
  \non
  \\
  && \rhot \to \rho
  \quad \hbox{weakly star in $\W{1,\infty}H\cap\H1V\cap\L\infty W$},
  \non
  \\
  && f_1'(\rhot) \to \xi
  \quad \hbox{weakly star in $\L\infty H$},
  \non
\Esist
at least for some subsequence $\tau_k\mathrel{\scriptstyle\searrow}0$.
Thanks to the Aubin-Lions lemma
(cf.~\cite[Thm.~5.1, p.~58]{Lions})
and to similar results to be found in \cite[Sect.~8, Cor.~4]{Simon},
we
{also}
deduce the following strong convergences:
\Bsist
  && \mut \to \mu
  \quad \hbox{strongly in $\C0H\cap\L2V$}
  \non
  \\
  && \rhot \to \rho
  \quad \hbox{strongly in $\C0V$}.
  \non
\Esist
In particular, {having recourse} to a \wk\ monotonicity technique
(see, e.g., \cite[Lemma~1.3, p.~42]{Barbu}),
we conclude that
$0<\rho<1$ and $\xi=f_1'(\rho)$ \aeQ.
The strong convergence 
{shown above}
also entails that
$f_2'(\rhot)$ converges to $f_2'(\rho)$, e.g., strongly in~$\C0H$
({because} $f'_2$ is Lipschitz continuous),
and that $\T\mut$ converges to~$\mu$, e.g., strongly in~$\L2H$.
Finally, a~combination of the above weak and strong convergence results
with the \holder\ and Sobolev inequalities yields that
\Bsist
  && \mut \, \dt\rhot \to \mu \, \dt\rho
  \quad \hbox{weakly in $\L1H$},
  \non
  \\
  && \rhot \, \dt\mut \to \rho \, \dt\mu
  \quad \hbox{weakly in $\L2{\Lx{3/2}}$}.
  \non
\Esist
Indeed, $\mut\to\mu$ strongly in $\Ldlq$,
$\dt\rhot\to\dt\rho$ weakly in $\Ldlq$,
$\rhot\to\rho$ strongly in $\C0{\Lx6}$,
and $\dt\mut\to\dt\mu$ weakly in~$\L2{\Lx2}$.
Therefore, it is \sfw\ to conclude that
the pair $(\mu,\rho)$ is a solution to problem~\Pbl\
having the desired regularity~\Regsoluz{, that is to say, Theorem~\ref{Esistenza} is proved}.

%%%%%%%%%%%%%%%%%%%%%%%%%%%%%%%%%%%%%%%%%%%%%%%%%%%%%%%%%%%%%%%%%%%%%%%%

\section{Uniqueness and boundedness}
\label{UniqBdd}
\setcounter{equation}{0}

In this section, we prove Theorem~\ref{Unicita} and Theorem~\ref{Limitatezza}.
We 
{first show our uniqueness result}.

\step Proof of Theorem~\ref{Unicita}

We take two solutions to problem \Pbl\ satisfying
\eqref{regsoluzbis} in addition to~\Regsoluz\ and 
label their components with the subscripts $1$ and~$2$; in the 
{following}, 
the values of constants $c$ 
{may depend on these solutions}.
Moreover, we choose constants such $M\geq0$ and $\rmin,\rmax\in(0,1)$ that
$\mu\leq M$ \aeQ\ and $\rmin\leq\rho_i\leq\rmax$ \aeQ, for $i=1,2$.
Finally, we {denote by} $L$ the Lipschitz constant of
the function $r\mapsto r-f'(r)$, $r\in[\rmin,\rmax]$.
{Having done this,} we write \eqref{prima} for both solutions
and take the difference.
Then, we set $\mu:=\mu_1-\mu_2$ and $\rho:=\rho_1-\rho_2$,
test the 
{resulting}
equality by~$\mu$,
and integrate, using the boundary condition.
Due to the identity:
\Beq
  \bigl\{
    (\dt\rho_1)\mu_1 + 2\rho_1 \dt\mu_1 
    - (\dt\rho_2)\mu_2 - 2\rho_2 \dt\mu_2 
  \bigr\} \, \mu
  = \dt(\rho_1 \mu^2) + 2 (\dt\mu_2) \, \rho \, \mu + \mu_2 \, (\dt\rho) \, \mu,
\Eeq
we obtain:
\Beq
  \iO \Bigl( \frac \eps 2 + \rho_1(t) \Bigr) \, |\mu(t)|^2
  + \intQt |\nabla\mu|^2
  = - 2 \intQt (\dt\mu_2) \, \rho \, \mu
  - \intQt \mu_2 \, (\dt\rho) \, \mu.
  \label{diffprima}
\Eeq
{Furthermore}, we write \eqref{seconda} for both solutions,
take the difference, test the 
{resulting}
equality by~$\dt\rho$,
and add $\rho\,\dt\rho$ to both sides, for convenience.
Then, we integrate, using the boundary condition, and easily obtain that
\Bsist
  && \delta \intQt |\dt\rho|^2
  + \frac 12 \, \normaV{\rho(t)}^2
  = \intQt \bigl( (\rho_1 - f'(\rho_1)) - (\rho_2 - f'(\rho_2)) + \mu \bigr) \, \dt\rho
  \non
  \\
  && \leq L \intQt |\rho| \, |\dt\rho|
  + \intQt |\mu| \, |\dt\rho| .
  \label{diffseconda}
\Esist

{Now, adding \eqref{diffprima} and \eqref{diffseconda} and taking into account that
$\rho_1$ is nonnegative}, 
we~get
\Bsist
  && \frac \eps 2 \iO |\mu(t)|^2
  + \intQt |\nabla\mu|^2
  + \delta \intQt |\dt\rho|^2
  + \frac 12 \, \normaV{\rho(t)}^2
  \non
  \\
  && \leq 2 \iot \normaH{\dt\mu_2(s)} \, \norma{\rho(s)}_{\Lq} \, \norma{\mu(s)}_{\Lq} \, ds
  + \iot \norma{\mu_2(s)}_{\Linfty} \, \normaH{\dt\rho(s)} \, \normaH{\mu(s)} \, ds
  \non
  \\
  && \quad {}
  + L \iot \normaH{\rho(s)} \, \normaH{\dt\rho(s)} \, ds
  + \iot \normaH{\mu(s)} \, \normaH{\dt\rho(s)} \, ds .
  \label{perunicita}
\Esist
{To estimate the first integral on the \rhs, we use the Sobolev inequality \eqref{sobolev} with $q=4$
 and $C$ the Sobolev constant, and we} 
{invoke}
the elementary Young inequality~\eqref{young} to
obtain that
\Bsist
  && \iot \normaH{\dt\mu_2(s)} \, \norma{\rho(s)}_{\Lq} \, \norma{\mu(s)}_{\Lq} \, ds
  \leq C^2 \iot \normaH{\dt\mu_2(s)} \, \normaV{\rho(s)} \, \normaV{\mu(s)} \, ds
  \non
  \\
  && \leq \frac 12 \iot \normaV{\mu(s)}^2 \, ds
  + \frac {C^4} 2 \iot \normaH{\dt\mu_2(s)}^2 \, \normaV{\rho(s)}^2 \, ds
  \non
  \\
  && = \frac 12 \intQt |\nabla\mu|^2 
  + c \intQt |\mu|^2
  + c \iot \normaH{\dt\mu_2(s)}^2 \, \normaV{\rho(s)}^2 \, ds .
  \non
\Esist
{The remainder of}
 the \rhs\ of \eqref{perunicita}
is estimated as follows:
\Bsist
  && \iot \norma{\mu_2(s)}_{\Linfty} \, \normaH{\dt\rho(s)} \, \normaH{\mu(s)} \, ds
  \non
  \\
  && \quad {}
  + L \iot \normaH{\rho(s)} \, \normaH{\dt\rho(s)} \, ds
  + \iot \normaH{\mu(s)} \, \normaH{\dt\rho(s)} \, ds
  \non
  \\
  && \leq \frac \delta 2 \intQt |\dt\rho|^2
  + c \iot \bigl( \normaH{\mu(s)}^2 + \normaV{\rho(s)}^2 \bigr) \, ds .
  \non
\Esist
Combining
{these}
estimates with \eqref{perunicita}, we immediately get
\Bsist
  && \frac \eps 2 \iO |\mu(t)|^2
  + \frac 12 \intQt |\nabla\mu|^2
  + \frac \delta 2 \intQt |\dt\rho|^2
  + \frac 12 \, \normaV{\rho(t)}^2
  \non
  \\
  && \leq c \iot \bigl( 1 + \normaH{\dt\mu_2(s)}^2 \bigr) \, \normaV{\rho(s)}^2 \, ds
  + c \iot \normaH{\mu(s)}^2 \, ds .
  \non
\Esist
Since the function $s\mapsto\normaH{\dt\mu_2(s)}^2$ belongs to $L^1(0,T)$,
we can apply the Gronwall lemma and deduce that both $\mu$ and $\rho$ vanish.
Hence, 
{the two solutions coincide.}\QED
%\marginpar{\footnotesizeQui ci starebbe bene il quadratino di fine dim.}

{We now turn to proving our boundedness result.}

\step Proof of Theorem~\ref{Limitatezza}

{Let $(\mu,\rho)$ be any solution to problem \Regsoluz\ and \Pbl\
whose initial data} {have, in addition to \Hpdati\,,  
the further properties}~\eqref{hpdatibis}.
{We show}
that the boundedness {claims in \eqref{regsoluzbis} actually hold true}.

{With a view to} proving that $\mu$ {satisfies the specified lower bound, we set:}
\Beq
  \muzs := \norma\muz_{\Linfty} = \mathop{\rm sup\,ess}_{x\in\Omega} \muz(x) \,;
  \label{supmuz}
\Eeq
we take any real constant $k$ such that $k\geq\muzs$;
and we introduce the auxiliary function $\chi_k\in {L^\infty (Q)}$ defined 
for a.a.\ $(x,t)\in Q$ by the formula:
\Beq
  \chik(x,t) = 1
  \quad \hbox{if $\mu(x,t)>k$},
  \aand
  \chik(x,t) = 0
  \quad \hbox{otherwise} .
  \non
\Eeq
Then, we test \eqref{prima} by $\mumkp$ and integrate over~$\Omega\times(0,t)$
for any $t\in(0,T)$. {The result is}:
\Bsist
  && \iO \Bigl(
    \frac \eps 2 + \rho(t)
  \Bigr) \, \bigl|\bigl( \mu(t) - k \bigr)^+\bigr|^2
  + \intQt |\nabla\mumkp|^2
  \non
  \\
  && = \intQt \dt\rho \, |\mumkp|^2
  - \intQt \dt\rho \, \mu \, \mumkp
  = - \intQt k \, \dt\rho \, \mumkp .
  \non
\Esist
{Given that $\rho$ is nonnegative, this equality}
 and the \holder\ inequality with ad hoc exponents lead~to:
\Bsist
  && \frac \eps 2 \, \normaH{(\mu(t)-k)^+}^2
  + \intQt |\nabla\mumkp|^2
  \non
  \\
  && \leq k \iot \norma{\chik(s)}_{\Lx{7/2}} \, \norma{\dt\rho(s)}_{\Lx{14/3}} \, \norma{\mumkp(s)}_{\Lx2} \, ds.
  \non
\Esist
Now, we use the Gronwall-Bellman lemma as in \cite[Lemma~A.4, p.~156]{Brezis},
and find that
\Bsist
  \hskip-1cm && \left\{
    \eps \norma\mumkp_{\C0H}^2
    + \intQ |\nabla\mumkp|^2
  \right\}^{1/2}
  \non
  \\
  \hskip-1cm && \leq {\frac{k}{\sqrt\eps}} \ioT {\norma{\chik(t)}_{\Lx{7/2}}} \, \norma{\dt\rho(t)}_{\Lx{14/3}} \, dt
 \non
  \\
  \hskip-1cm &&  \leq {\frac{k}{\sqrt\eps}}  \, \norma{\dt\rho}_{\L{7/3}{\Lx{14/3}}} \,
    \norma\chik_{\L{7/4}{\Lx{7/2}}} \,.
  \non
\Esist
Next, we observe that the interpolation inequality~\eqref{Interpolazione}
(with $p=2$, $q=6$, $r=14/3$, and $\theta=1/7$),
together with the Sobolev inequality \eqref{sobolev}, gives that
\Bsist
  && \norma v_{\L{7/3}{\Lx{14/3}}}
  \leq \Bigl( \ioT \norma{v(t)}_{\Lx2}^{1/3} \, \norma{v(t)}_{\Lx6}^2 \, dt \Bigr)^{3/7}
  \non
  \\
  && = \norma v_{\L\infty H}^{1/7} \Bigl( \ioT \norma{v(t)}_{\Lx6}^2 \Bigr)^{3/7}
  = \norma v_{\L\infty H}^{1/7} \, \norma v_{\L2{\Lx6}}^{6/7}
  \non
  \\
  && \leq c \norma v_{\L\infty H}^{1/7} \, \norma v_{\L2V}^{6/7}, 
  \non
\Esist
{and we denote by $\Dzero$ the rightmost side of this inequality chain, when evaluated for}
$v=\dt\rho$.
We also remark that
\Bsist
  && \norma\chik_{\L{7/4}{\Lx{7/2}}}
  = \Bigl\{
    \ioT \Bigl(
      \iO |\chik(x,t)|^{7/2} \, dx
    \Bigr) ^{1/2} dt
  \Bigr\}^{4/7}
  \non
  \\
  && = \Bigl\{
    \ioT \Bigl(
      \iO |\chik(x,t)|^4 \, dx
    \Bigr) ^{1/2} dt
  \Bigr\}^{\frac12\cdot\frac87}
  = \norma\chik_{\Ldlq}^{8/7} .
  \non
\Esist
Hence, our estimate for $\mumkp$ yields the following basic inequality:
\Beq
  \Norma\mumkp \leq k \Duno \norma\chik_{\Ldlq}^{8/7}
  \quad \hbox{for every $k\geq\muzs$},
  \label{bravopier}
\Eeq
where $\Duno=\Dzero/\min\graffe{\eps,1}$, and where
the norm $\Norma{\cdot}$ is defined by
\Beq
  \Norma v^2 := \sup_{t\in[0,T]}\normaH{v{(t)}}^2 + \int_Q |\nabla v|^2
  \quad \hbox{for $v\in\C0H\cap\L2V$}.
  \non
\Eeq
We notice that the Sobolev inequality \eqref{sobolev} implies that 
\Beq
  \norma v_{\Ldlq} \leq \Ddue \Norma v
  \quad \hbox{for every $v\in\C0H\cap\L2V$},
  \label{dasobolev}
\Eeq
where $\Ddue$ depends on $\Omega$ and~$T$, only.
At this point, we {select a sequence $\graffe{k_j}$
depending on a real parameter $m>1$ as follows:}
\Beq
%  M := m \muzs 
%  \aand
  k_j := M \bigl( 2 - 2^{-j} \bigr)
  \quad \hbox{for $j=0,1,\dots$,}\quad\textrm{with}\quad   M := m \muzs ;
  \label{sceltakj}
\Eeq
note that $k_0=M>\muzs$.
Then, owing to \eqref{bravopier} and~\eqref{dasobolev}, it is not difficult to check that
\Bsist
  && \bigl( k_{j+1} - k_j \bigr) \, \norma{\chi_{k_{j+1}}}_{\Ldlq}
  \leq \norma{(\mu-k_j)^+}_{\Ldlq}
  \leq \Ddue \Norma{(\mu-k_j)^+}
  \non
  \\
  && \leq k_j \, \Duno \Ddue \norma{\chi_{k_j}}_{\Ldlq}^{8/7} .
  \label{ricorrente}
\Esist
Therefore, if we set
\Beq
  S_j := \norma{\chi_{k_j}}_{\Ldlq}
  \quad \hbox{for $j=0,1,\dots$,}
  \non
\Eeq
then  {the following inequality holds:}
\Beq
  S_{j+1} 
  \leq \frac {k_j}{k_{j+1}-k_j} \, \Duno\Ddue S_j^{8/7}
  \leq 4 \Duno \Ddue \, 2^j S_j^{8/7}
  \quad \hbox{for $j=0,1,\dots$}.
  \non
\Eeq
Using \cite[Lemma 5.6, p.~95]{LSU},
we conclude that $S_j\to0$ as $j\to\infty$, provided that
\Beq
  S_0 = \norma{\chi_{k_0}}_{\Ldlq}
  \leq (4\Duno\Ddue)^{-7} 2^{-49}.
  \label{daLSU}
\Eeq
On the other hand, we notice that $\chi_{k_0}=\chi_M$,
and we recall that $M>\muzs$ and $m=M/\muzs$, by \eqref{sceltakj}.
Moreover, we observe that 
$\chi_M=1<(\mu-\muzs)/(M-\muzs)$ when $\mu>M$,
and that $\chi_M=0$ otherwise.
Therefore, using~\eqref{dasobolev} and \eqref{bravopier} with $k=k_0=M$, we have:
\Bsist
  && S_0 
  \leq \frac 1 {M-\muzs} \, \norma{(\mu-\muzs)^+}_{\Ldlq}
  \leq \frac \Ddue {M-\muzs} \, \Norma{(\mu-\muzs)^+}
  \non
  \\
  && \frac {\Duno\Ddue} {m-1} \, \norma{\chi_{\muzs}}_{\Ldlq}^{8/7}
  \leq \frac {\Duno\Ddue} {m-1} \, |\Omega|^{\frac14\cdot\frac87} \, T^{\frac12\cdot\frac87}.
  \non
\Esist
{We are now in a position to} choose
$m:=1+\Duno\Ddue|\Omega|^{2/7}T^{4/7}(4\Duno\Ddue)^72^{49}$.
Then, $m>1$ and \eqref{daLSU} is satisfied.
Consequently,
\Beq
  \norma{\chi_{2M}}_{\Ldlq} = \lim_{j\to\infty} S_j = 0,
  \non
\Eeq
due to {Beppo Levi's Monotone Convergence Theorem}.
This implies 
{that $\mu\leq2M$ \aeQ\,,}
 and the boundedness of $\mu$ 
claimed in \eqref{regsoluzbis} is established.

{We are left with the task of proving that the limitations for $\rho$ in~\eqref{regsoluzbis} do hold}.
{We find it convenient to set: $\rhomin:=\inf_\Omega\rhoz$
(recall that we assumed $\rhomin$ to be strictly positive)}.
Moreover, we rewrite \eqref{seconda} in the~form:
\Beq
  \delta \dt\rho - \Delta\rho + f_1'(\rho) = g,
  \quad \hbox{where} \quad
  g := \mu - f_2'(\rho),
  \label{riseconda}
\Eeq
and we notice that $g\in\LQ\infty$, in view of the above proof and~\eqref{hpfi}.
{Consequently, in view also of}~\eqref{hpfp},
we can choose $\rmin\in(0,\rhomin)$ such that $f_1'(\rmin)\leq g$ \aeQ.
Then, we test \eqref{riseconda} by $-\rhomrm$ and deduce~that
\Bsist
  && \frac \delta 2 \iO |\rhomrm(t)|^2
  + \intQt |\nabla\rhomrm|^2
  - \intQt \bigl( f_1'(\rho) - f_1'(\rmin) \bigl) \rhomrm
  \non
  \\
  && = \intQt \bigl( f_1'(\rmin) - g \bigl) \rhomrm
  \leq 0 .
  \non
\Esist
We conclude that $\rhomrm=0$ and $\rho\geq\rmin$ \aeQ.
In a similar way, for a suitable $\rmax<1$,
we show that $\rho\leq\rmax$ \aeQ\
by testing \eqref{riseconda} by $(\rho-1+\rmax)^+$.
 {We conclude that solutions satisfy all of the requirements stated in \eqref{regsoluzbis}}.\QED
%\marginpar{\footnotesize Anche qui ci metterei il quadratino di fine dim.}

%%%%%%%%%%%%%%%%%%%%%%%%%%%%%%%%%%%%%%%%%%%%%%%%%%%%%%%%%%%%%%%%%%%%%%%%

\section{\Loti\ \bhv}
\label{LongTime}
\setcounter{equation}{0}

In this section, we prove Theorem~\ref{Longtime}.
{To this end,}
 we fix any solution $(\mu,\rho)$ to problem \Pbl.
Our proof of the properties of the $\omega$-limit $\omega(\mu,\rho)$
relies on a number of a~priori estimates for $(\mu,\rho)$, and on a \wk\ tool.
{For} $(\muo,\rhoo)$ any element of $\omega(\mu,\rho)$,
and $\graffe{t_n}$ a corresponding {time sequence of type}~\eqref{defomegalim},
we~set
\Beq
  \mun(t) := \mu(t_n+t),\quad
  \rhon(t) := \rho(t_n+t)
  \quad \hbox{for $t\geq0$},
  \label{munrhon}
\Eeq
and we study the sequence $\graffe{(\mun,\rhon)}$ 
on a fixed finite time interval~$[0,T]$.
Clearly, the pair $(\mun,\rhon)$
enjoys the same regularity 
{as $(\mu,\rho)$, and solves}
 the equations
\Bsist
  & (\eps + 2\rhon) \dt\mun + \mun \, \dt\rhon - \Delta\mun = 0
  & \quad \aeQ
  \label{primalongn}
  \\
  & \delta \dt\rhon - \Delta\rhon + f'(\rhon) = \mun
  & \quad \aeQ ;
  \label{secondalongn}  
\Esist
moreover, it satisfies the homogeneous Neumann boundary conditions
and the Cauchy conditions:
\Beq
  \mun(0) = \mu(t_n)
  \aand
  \rhon(0) = \rho(t_n)
  \quad \aeO .
  \label{cauchylong}
\Eeq
Our argument also relies on 
{two basic identities, to be proved}
in the next lemma.

\Blem
\label{Juergen}
The following identities hold:
\Bsist
  && \delta (\dt\rho)^2 - \dt\rho \, \Delta\rho + f'(\rho) \, \dt\rho
  = \eps \dt\mu + 2 \dt(\rho\mu) - \Delta\mu
  \vphantom\int\quad \aeQ ;
  \label{juergen}
  \\
  && \delta \iot \normaH{\dt\rho(s)}^2 \, ds
  + \frac 12 \, \normaH{\nabla\rho(t)}^2
  + \iO f(\rho(t))
  \non
  \\
  && = \frac 12 \, \normaH{\nabla\rhoz}^2
  + \iO f(\rhoz)
  + \eps \iO \mu(t)
  - \eps \iO \muz 
  + 2 \iO (\rho\mu)(t) 
  - 2 \iO \rhoz\muz,
  \qquad
  \label{intjuergen}
\Esist
for every $t\in[0,T]$.
\Elem

\par\noindent{\bf Proof.} 
We have from~\eqref{prima} that
\Beq
  \mu \, \dt\rho
  = \Delta\mu - \eps \dt\mu - 2\rho \, \dt\mu
  = \Delta\mu - \eps \dt\mu - 2\dt(\rho\mu) + 2 \dt\rho \, \mu \,.
  \non
\Eeq
{By a simple rearrangement,}
we deduce that
\Beq
  \mu \, \dt\rho
  = \eps \dt\mu + 2 \dt(\rho\mu) - \Delta\mu .
  \non
\Eeq
On the other hand, 
{multiplication of}
\eqref{seconda} by $\dt\rho$ yields:
\Beq
  \delta (\dt\rho)^2 - \dt\rho \, \Delta\rho + f'(\rho) \, \dt\rho
  = \mu \, \dt\rho, 
  \non
\Eeq
so that \eqref{juergen} immediately follows by comparison.
Next, identity~\eqref{intjuergen} is arrived at by integrating \eqref{juergen} over $\Omega\times(0,t)$
and by noting 
{that, in view of the homogeneous Neumann boundary condition,
 $\Delta\mu$ does not contribute to the integral.}\QED

{We proceed with proving} some a~priori estimates. {In so doing, we depart from our general rule, and write $c$ for constants that do not depend
on the final time~$T$, although they are allowed to depend on the element
$(\muo,\rhoo)$ of the $\omega$-limit 
under consideration. Whenever a dependence of $c$ on the parameter~$T$ cannot be excluded, we stress this possibility by writing $c_T$.}
{Morever, without any loss of generality, we
assume that $\eps\leq1$.}

\step First a priori estimate

Just as we did for \eqref{perprimastima}, we immediately deduce~that
\Beq
  \iot \normaH{\nabla\mu(s)}^2 \, ds
  + {\frac\eps 2} \normaH{\mu(t)}^2 
  + \iO (\rho\mu^2)(t)
  \leq c
  \quad \hbox{for every $t>0$}.
  \label{perprimalong}
\Eeq
This implies, in particular, that
\Beq
  \int_0^{+\infty} \normaH{\nabla\mu(t)}^2 \, dt < +\infty
  \aand
  \norma\mun_{\L\infty H} \leq c.
  \label{primalong}
\Eeq

\step Second a priori estimate

We recall \eqref{intjuergen} and estimate some terms of its \rhs.
We have:
\Beq
  \eps \iO \mu(t)
  + 2 \iO (\rho\mu)(t) 
  \leq
  \eps^{1/2} \iO \mu(t)
  + 2 \iO (\rho^{1/2}\mu)(t)
  \leq 2 |\Omega| + \eps \normaH{\mu(t)}^2 
  + \iO (\rho\mu^2)(t) .
  \non
\Eeq
On the other hand, \eqref{perprimalong} holds
and $f$ is bounded from below.
Hence, \eqref{intjuergen} yields:
\Beq
  \delta \iot \normaH{\dt\rho(s)}^2 \, ds
  + \frac 12 \, \normaH{\nabla\rho(t)}^2
  \leq c 
  \quad \hbox{for every $t>0$},
  \label{persecondalong}
\Eeq
whence we have that
\Beq
  \int_0^{+\infty} \normaH{\dt\rho(t)}^2 \, dt < +\infty
  \aand
  \norma\rhon_{\L\infty V} \leq c.
  \label{secondalong}
\Eeq

\step Third a priori estimate

We formally test \eqref{seconda} by $-\Delta\rho$
and integrate over $\Omega\times(t_n,t_n+t)$.
Due to the convexity of $f_1$ and the boundedness of~$f_2''$
{(compare with the derivation of~\eqref{dadtsecondanl}),}
we~get
\Bsist
  && \frac \delta 2 \, \normaH{\nabla\rho(t_n+t)}^2
  + \intQtnt |\Delta\rho|^2
  \non
  \\
  && \leq \frac \delta 2 \, \normaH{\nabla\rho(t_n)}^2
  + c \intQtnt |\nabla\rho|^2 
  + \frac 12 \intQtnt |\Delta\rho|^2
  + \frac 12 \intQtnt |\mu|^2 .
  \non
\Esist
We note that the first term on the \rhs\ is bounded,
since $\rho(t_n)$ is weakly convergent to 
$\rhoo$ in~$V$.
Hence, owing to \eqref{persecondalong} and~\eqref{perprimalong},
we can conclude that
\Beq
  \normaH{\nabla\rho(t_n+t)}^2
  + \intQtnt |\Delta\rho|^2
  \leq c_T
  \quad \hbox{for every $t\in[0,T]$ and every $n$}.
  \non
\Eeq
By comparison with~\eqref{seconda}, and by exploiting elliptic regularity, we deduce~that
\Beq
  \itnt \normaH{f_1'(\rho)}^2 \, dt
  + \itnt \normaW{\rho(t)}^2 \, dt \leq c_T \,.
  \non
\Eeq
In terms of $\rhon$, all this {reads:}
\Beq
  \norma\rhon_{\L2W} \leq c_T
  \aand
  \norma{f_1'(\rhon)}_{\L2H} \leq c_T \,.
  \label{terzalong}
\Eeq

\step Fourth a priori estimate

By rewriting \eqref{primalongn} in the form
\Beq
  \dt\mun
  = - \frac {\mun \, \dt\rhon} {\eps + 2\rhon}
  + \frac {\Delta\mun} {\eps + 2\rhon}, 
  \label{reprimalong}
\Eeq
and owing to the homogeneous Neumann boundary condition,
we {are entitled to write the following} equation in~$\Vp$
in~the framework of the Hilbert triplet $(V,H,\Vp)$:
\Beq
 \iO \dt\mu(t) \, v
  = - \iO \frac {\mun(t) \, \dt\rhon(t) \, v} {\eps + 2\rhon(t)}
  - \iO \nabla\mun(t) \cdot \nabla \frac v {\eps + 2\rhon(t)}\,,
  \label{testata}
\Eeq
\aat\ and for every $v\in V$.
Starting from this equation, 
we can prove a bound for $\dt\mun$ in $\L p\Vp$ for some $p>1$.
For a while, we argue \aat\ and estimate 
each term on the \rhs\ separately.
Owing to the Sobolev inequality~\eqref{sobolev},
we~get
\Bsist
  && \Bigl| \iO \frac {\mun(t) \, \dt\rhon(t) \, v} {\eps + 2\rhon(t)} \Bigr|
  \leq \eps^{-1} \norma{\mun(t)}_{\Lx4} \, \norma{\dt\rho(t)}_{\Lx2} \, \norma v_{\Lx4}
  \non
  \\
  && \leq c \norma{\mun(t)}_{\Lx4} \, \norma{\dt\rho(t)}_{\Lx2} \, \normaV v .
  \label{pernormaVpa}
\Esist
On the other hand, using the interpolation inequality \eqref{Interpolazione} 
and the Sobolev inequality once more, we obtain that
\Beq
  \norma{\mun(t)}_{\Lx4}
  \leq \norma{\mun(t)}_{\Lx6}^{3/4} \, \norma{\mun(t)}_{\Lx2}^{1/4}
  \leq c \normaV{\mun(t)}^{3/4} \, \normaH{\mun(t)}^{1/4}.
  \non
\Eeq
Consequently, by accounting for the $L^\infty$ bound of~\eqref{primalong},
we derive from~\eqref{pernormaVpa} that
\Beq
  \Bigl| \iO \frac {\mun(t) \, \dt\rhon(t) \, v} {\eps + 2\rhon(t)} \Bigr|
  \leq c \normaV{\mun(t)}^{3/4} \, \normaH{\dt\rho(t)} \, \normaV v.
  \label{normaVpa}
\Eeq
{This takes care of the first addendum in on the \rhs\ of~\eqref{testata}. As to the second, we} have:
\Beq
  \Bigl| \iO \nabla\mun(t) \cdot \nabla \frac v {\eps + 2\rhon(t)} \Bigr|
  \leq \normaV{\mun(t)} \, \normaH{\nabla\bigl((\eps+2\rhon(t))^{-1}\,v\bigr)},
  \label{pernormaVpb}
\Eeq
where 
{it remains for us}
to estimate the last norm.
{By applying the Leibniz rule and making use of
the \holder\  interpolation and the Sobolev inequalities}, we~get:
\Bsist
  && \normaH{\nabla((\eps+2\rhon(t))^{-1}\,v)}
  \leq 2\eps^{-2} \normaH{v\nabla\rhon(t)}
  + \eps^{-1} \normaH{\nabla v}
  \non
  \\
  && \leq c \bigl( \norma v_{\Lx4} \, \norma{\nabla\rhon(t)}_{\Lx4} + \normaH{\nabla v} \bigr)
  \leq c \bigl( \norma{\nabla\rhon(t)}_{\Lx4} + 1 \bigr) \normaV v
  \non
  \\
  && \leq c \bigl( \norma{\nabla\rhon(t)}_{\Lx6}^{3/4} \, \norma{\nabla\rhon(t)}_{\Lx2}^{1/4} + 1 \bigr) \normaV v
  \leq c \bigl( \normaW{\rhon(t)}^{3/4} \, \normaV{\rhon(t)}^{1/4} + 1 \bigr) \normaV v \,.
  \non
\Esist
Hence, on accounting for the $L^\infty$ bound of~\eqref{secondalong}, we see that \eqref{pernormaVpb} becomes:
\Beq
  \Bigl| \iO \nabla\mun(t) \cdot \nabla \frac v {\eps + 2\rhon(t)} \Bigr|
  \leq c \normaV{\mun(t)} \bigl( \normaW{\rhon(t)}^{3/4} + 1 \bigr) \, \normaV v \,.
  \label{normaVpb}
\Eeq
Since $v\in V$ is arbitrary,
by combining \eqref{testata}, \eqref{normaVpa}, and~\eqref{normaVpb},
we arrive at:
\Beq
  \normaVp{\dt\mun(t)}
  \leq c \bigl(
    \normaV{\mun(t)}^{3/4} \, \normaH{\dt\rho(t)}
    + \normaV{\mun(t)} \, \normaW{\rhon(t)}^{3/4} + \normaV{\mun(t)}
  \bigr),
  \quad
  \label{normaVp}
\Eeq
\aat\ and for every~$n$; {by estimating each term on the \rhs,
we are going to} find some $p>1$  such that
\Beq
  \norma{\dt\mun}_{\L p\Vp} \leq c_T.
  \label{quartalong}
\Eeq
Now, due to~\eqref{primalong}, the last term is bounded in~$L^2(0,T)$. {As to the first,} we observe that the functions
\Beq
  t \mapsto \normaV{\mun(t)}^{3/4} 
  \aand
  t \mapsto \normaH{\dt\rho(t)}
  \non
\Eeq
are bounded, respectively,  in $L^{8/3}(0,T)$ by \eqref{primalong}, and in~$L^2(0,T)$ 
 by~\eqref{secondalong};
hence, their product is bounded in $L^{8/7}(0,T)$, by the \holder\ inequality.
{Finally, the middle term of~\eqref{normaVp} can be treated in a similar way, with the use of
\eqref{primalong} and the first inequality in~\eqref{terzalong}.}
Therefore,
\eqref{quartalong} {does hold,} with $p=8/7$.

\step Conclusion

Our estimates \eqref{primalong}, \eqref{secondalong}, and~\eqref{quartalong}
ensure that $(\mu,\rho)$ is a bounded and weakly continuous
$(H,V)$-valued function.
Hence, the first part of the statement {in Theorem~\ref{Longtime}} follows from the general theory
(see, e.g., \cite[p.~12]{Haraux}).
We pass to the study of the $\omega$-limit.

Recalling \eqref{terzalong} 
and using standard weak and weak star compactness results, we see that there is
a triplet $(\mui,\rhoi,\phii)$ such that
\Bsist
  \hskip-1cm & \mun \to \mui
  & \quad \hbox{weakly {star} in $\L\infty H\cap\L2V$},
  \label{convmu}
  \\
  \hskip-1cm & \rhon \to \rhoi
  & \quad \hbox{weakly star in $\H1H {\cap \L\infty V} \cap\L2W$},
  \label{convrho}
  \\
  \hskip-1cm & f_1'(\rhon) \to \phii
  & \quad \hbox{weakly in $\L2H$},
  \label{convf}
\Esist
at least for some subsequence.
Our first aim is to prove that $\mui$ is a 
{nonnegative}
constant,
i.e., that $\mui(x,t)=\mus$ \aaQ\ for some $\mus\in[0,+\infty)${; and, to prove
that} $\rhoi$ is time independent, i.e., that $\rhoi(t)=\rhos$ \aat\
for some $\rhos\in W$.
Secondly, we want to prove that the pair $(\mus,\rhos)$ found in such a way
is indeed a steady state and coincides with the given pair~$(\muo,\rhoo)$.

From the first bounds of \eqref{primalong} and~\eqref{secondalong}, we immediately deduce~that
\Beq
  |\nabla\mun| \to 0
  \aand
  \dt\rhon \to 0
  \quad \hbox{strongly in $\L2H$}.
  \non
\Eeq
This implies that $\mui$ is space independent
and $\rhoi$ is time independent.
Thus, we can write $\rhoi(t)=\rhos$ \aat,
for some $\rhos\in W$.
Moreover, \eqref{convrho}~implies strong convergence:
\Beq
  \rhon \to \rhoi 
  \quad \hbox{strongly in $\ {\C0H \cap}\L2V$}
  \label{strongconvrho}
\Eeq
{(see, e.g., \cite[Sect.~8, Cor.~4]{Simon})}.
Therefore, $f'_2(\rhon)$ converges to $f'_2(\rhoi)$,
e.g., strongly in $\L2H$, and it is clear~that
\Beq
  - \Delta\rhoi + \phii = \mui - f'_2(\rhoi)
  \quad \aeQ.
  \label{limsecondalong}
\Eeq
Collecting \eqref{strongconvrho} and~\eqref{convf}, 
and recalling 
{Lemma~1.3, p.~42, in
\cite{Barbu},}
we conclude that
$0<\rhoi<1$ and $\phii=f_1'(\rhoi)$ \aeQ.
Therefore, \eqref{limsecondalong}~becomes:
\Beq
  0 < \rhos < 1
  \aand
  - \Delta\rhos + f'(\rhos) = \mui
  \quad \aeQ,
  \non
\Eeq
and we deduce that $\mui$ is time independent as well.
Thus, $\mui(x,t)=\mus$ \aaQ\ for some constant~$\mus$.
Furthermore, $\mus$~is 
{nonnegative, since $\mun\geq0$ for every~$n$.}
This concludes the proof that $(\mus,\rhos)$ is a steady state.

It remains { for us} to show that $(\mus,\rhos)$ coincides with~$(\muo,\rhoo)$.
From~\eqref{strongconvrho} we see that
$\rhon(0)$ converges strongly in~${H}$ to~$\rhoi(0)=\rhos$;
on the other hand, $\rhon(0)=\rho(t_n)$ converges weakly in $V$ to $\rhoo$,
 by assumption; hence, $\rhos=\rhoo$.
A~similar argument holds for $\mus$ and $\muo$, {because $\mu_n$ converges strongly} in $\C0\Vp$.
Indeed, $\mun$ is bounded in $\L\infty H$, by~\eqref{primalong};
on the other hand, \eqref{quartalong}~holds with $p=8/7>1$, 
and the embedding $H\subset\Vp$ is compact; hence, the desired convergence follows from~\cite[Sect.~8, Cor.~4]{Simon}.
This completes the proof {of Theorem~\ref{Longtime}}.

%%%%%%%%%%%%%%%%%%%%%%%%%%%%%%%%%%%%%%%%%%%%%%%%%%%%%%%%%%%%%%%%%%%%%%%%

%%%%%%%%%%%%%%%%%%%%%%%%%%%%%%%%%
%% bibliography
%%%%%%%%%%%%%%%%%%%%%%%%%%%%%%%%%

\vspace{3truemm}

\Begin{thebibliography}{10}

\bibitem{Barbu}
V. Barbu,
``Nonlinear semigroups and differential equations in Banach spaces'',
Noord\-hoff,
Leyden,
1976.

\bibitem{bds}
{E. Bonetti, W. Dreyer, G. Schimperna, Global solutions to 
a generalized {C}ahn-{H}illiard equation with viscosity,
{\it Adv. Differential Equations} {\bf 8} (2003) 231-256.}

\bibitem{Brezis}
H. Brezis,
``Op\'erateurs maximaux monotones et semi-groupes de contractions
dans les espaces de Hilbert'',
North-Holland Math. Stud.
{\bf 5},
North-Holland,
Amsterdam,
1973.

{\bibitem{CN} B.D. Coleman, W. Noll, 
The thermodynamics of elastic materials with heat conduction and viscosity, 
{\it Arch. Rational Mech. Anal.} {\bf 13} (1963) 167-178.} 

\bibitem{CGPS} 
P. Colli, G. Gilardi, P. Podio-Guidugli, J. Sprekels,
Existence and uniqueness of a global-in-time solution
to a phase segregation problem of the Allen-Cahn type, 
{\it Math. Models Methods Appl. Sci.} {\bf 20} (2010)
519-541.

\bibitem{CGPS2} 
P. Colli, G. Gilardi, P. Podio-Guidugli, J. Sprekels,
A temperature-dependent phase segregation problem 
of the Allen-Cahn type, {\it Adv. Math. Sci. Appl.} 
{\bf 20} (2010) 
{219-234.}

\bibitem{Fremond}
M. Fr\'emond,
``Non-smooth Thermomechanics'',
Springer-Verlag, Berlin, 2002.

{\bibitem{FG} 
E. Fried, M.E. Gurtin, 
Continuum theory of thermally induced phase transitions based on an order 
parameter, {\it Phys. D} {\bf 68} (1993) 326-343.}

\bibitem{GR} 
G. Gilardi, E. Rocca,
Well posedness and long time behaviour for a singular phase field system of conserved type,
{\it IMA J. Appl. Math.} {\bf 72} (2007) 498-530. 

\bibitem{Gurtin} 
M. Gurtin, Generalized Ginzburg-Landau and
Cahn-Hilliard equations based on a microforce balance,
{\it Phys.~D\/} {\bf 92} (1996) 178-192.

\bibitem{Haraux}
A. Haraux,
``Syst\`emes Dynamiques Dissipatifs et Applications'',
RMA Res. Notes Appl. Math.,
{\bf 17},
Masson,
Paris,
1991.

\bibitem{LSU}
O.A. Lady\v zenskaja, V.A. Solonnikov, and N.N. Ural'ceva:
``Linear and quasilinear equations of parabolic type'',
Trans. Amer. Math. Soc., {\bf 23},
Amer. Math. Soc., Providence, RI,
1968.

\bibitem{Lions}
J.L. Lions,
``Quelques m\'ethodes de r\'esolution des probl\`emes aux limites non
lin\'eaires'',
Dunod Gauthier--Villars,
Paris,
1969.

\bibitem{mr} 
{A. Miranville, A Rougirel, Local and asymptotic analysis of the flow generated by the
{C}ahn-{H}illiard-{G}urtin equations,  {\it Z. Angew. Math. Phys.} {\bf 57} (2006) 244-268.}

\bibitem{Podio}
P. Podio-Guidugli, 
Models of phase segregation and diffusion of atomic species on a lattice,
{\it Ric. Mat.} {\bf 55} (2006) 105-118.

\bibitem{Ros} 
{R. Rossi, On two classes of generalized viscous {C}ahn-{H}illiard
equations, {\it Commun. Pure Appl. Anal.} {\bf 4} (2005) 405-430.}

\bibitem{Simon}
J. Simon,
{Compact sets in the space $L^p(0,T; B)$},
{\it Ann. Mat. Pura Appl.} {\bf 146} (1987) 65--96.

\End{thebibliography}

\End{document}

\bye